\documentclass[12pt]{article}
\usepackage{latexsym, amsfonts, amscd, amssymb, verbatim, amsmath,enumerate}

\begin{document}

\title{Toward Homological Characterization of Semirings by e-Injective
Semimodules}
\author{J. Y. Abuhlail\thanks{%
Department of Mathematics and Statistics, King Fahd University of Petroleum
and Minerals, Dhahran, Saudi Arabia. Email: \texttt{abuhlail@kfupm.edu.sa}},
S. N. Il'in\thanks{%
Lobachevsky Institute of Mathematics and Mechanics, Kazan (Volga Region) Federal
University, Kazan, Tatarstan, Russia. Email: \texttt{sergey.ilyin@kpfu.ru}},
Y. Katsov\thanks{%
Department of Mathematics, Hanover College, Hanover, IN 47243--0890, USA.
Email: \texttt{katsov@hanover.edu}} and T.\thinspace G.~Nam\thanks{%
Institute of Mathematics, VAST, 18 Hoang Quoc Viet, Cau Giay, Hanoi,
Vietnam. Email: \texttt{tgnam@math.ac.vn}}}
\date{\vspace{-12mm}}
\maketitle

\begin{abstract}
In this paper, we introduce and study e-injective semimodules, in particular
over additively idempotent semirings. We completely characterize semirings
all of whose semimodules are e-injective, describe semirings all of whose
projective semimodules are e-injective, and characterize one-sided
Noetherian rings in terms of direct sums of e-injective semimodules. Also,
we give complete characterizations of bounded distributive lattices,
subtractive semirings, and simple semirings, all of whose cyclic (finitely
generated) semimodules are e-injective.


\end{abstract}

\section{Introduction}

Semirings and semimodules, and their applications, arise in various branches
of Mathematics, Computer Science, Physics, as well as many other areas of
modern science (see, for instance, \cite{golan:sata} and \cite{glaz:agtl}).
In the recent years, there has been a substantial amount of interest in
additively idempotent semirings --- among which the Boolean semifield,
tropical semifields, and coordinate semirings of tropical varieties
represent a set of well-known examples --- originated in several extremely
interesting, \textquotedblleft nontraditional\textquotedblright\ contexts as
Tropical Geometry \cite{rst:fsitg}, Tropical Algebra \cite{ir:sa}, $\mathbb{F%
}_{1}$-Geometry \cite{cc:sofazf}, and the Geometry of Blueprints \cite%
{lor:tgob}, for example. Also, in the last decade, motivated by the Riemann
Hypothesis \cite{cc:sofazf} and tropical varieties \cite{rst:fsitg} and \cite%
{gg:eotv}, several mathematicians have studied from different points of view
semiring schemes, in particular, in \cite{jun:ccoss}, sheaves and
homological methods on semiring schemes have been considered.

In the same time, homological characterization/classification of rings by
properties of suitable classes (categories) of modules over them ( see,
\textit{e.g.}, \cite{lam:lomar}) constitutes one of the most sustained
interests and important achievements of Homological Algebra. Inspired by
this, during the last three decades, a good number of important results
related to this genre have been obtained in different non-abelian settings
as, for example, in the homological classification of monoids \cite%
{kilp-kn-mik:maac} and distributive lattices \cite{fof:podl}. As algebraic
objects, semirings are certainly the most natural generalization of such (at
first glance different) algebraic systems as rings and bounded distributive
lattices, and therefore, they form an extremely interesting, natural, and
important, non-abelian/non-additive setting for furthering of the
homological characterization, \textit{i.e.}, characterizing semirings by
properties of suitable classes (categories) of semimodules over them. In
fact, this is an ongoing project of an substantial interest (see, \textit{%
e.g.}, \cite{kat:tpaieosoars}, \cite{kat:thcos}, \cite{kat:ofsos}, \cite%
{il'in:svwasai}, \cite{knt:ossss}, \cite{il'in:dsoisadpopsos}, \cite%
{kn:meahcos}, \cite{il'in:v-s}, \cite{ikn:thstososaowcsap}, and \cite%
{aikn:ovsasaowcsai}). In all studies regarding the homological
characterization, the concepts of `injectivity' and `projectivity' of
objects --- $R$-modules, $S$-acts, $S$-semimodules, etc.--- play the most
leading role. Perhaps one of the most important achievements of homological
algebra are due to the fact that for abelian categories of modules $_{R}%
\mathcal{M}$ over a ring $R$\ the both --- universal algebra and homological
algebra --- approaches lead to the identical classes of modules (see, for
example, \cite[Sections 1.2A and 1.3A]{lam:lomar} or \cite[Sect. II.6.9]%
{gelman:moha}). However, as we will see later on, in generally additive, but
non-abelian, setting of categories of semimodules $_{S}\mathcal{M}$ over a
semiring $S$, these two approaches lead to two different classes of
semimodules. In non-abelian (even in non-additive) settings, there have been
obtained a good number of quite interesting and important results connected
with the concepts of `projectivity' and `injectivity' based on the universal
algebra approach (see, \textit{e.g.}, \cite{kilp-kn-mik:maac}). In contrast
to this, in the present paper, we initiate investigations related to
injective and projective semimodules defined by using the second approach,
\textit{i.e.}, from the homological algebra point of view, heavily based on
the fundamental concepts of `extensions' and `short exact sequences' of
modules, that, in turn, lead us to the concepts of `e-injectivity' and
`e-projectivity' of semimodules. We should mention that the latter concepts
have been earlier somewise considered by some authors in semimodule settings
under different terminology (see, for example, \cite{Tak:eos}, \cite{Tak81},
\cite{Tak1982a}, \cite{Tak82}, \cite{p:eosatfe}, \cite{Abu2014-a}, and \cite%
{Abu2014-b}), but our approach, nevertheless, is slightly different and, we
hope, better reflects the \textquotedblleft homological\textquotedblright\
spirit\ of the matter and based on it proofs of the\ obtained results.

The paper is organized as follows. In Section 2, for the reader's
convenience, we provide all subsequently necessary notions and facts on
semirings and semimodules.

In Section 3, we introduce the concepts of e-injectivity (-projectivity) of
semimodules and establish some fundamental facts regarding them we use in a
sequence. Among results of this section, we single out, in our view quite
interesting and useful, Proposition 3.4 and Corollary 3.5 that provide us
with a very useful tool to construct new (e-)injective semimodules from
known ones.

In Section 4, we completely characterize e-injective semimodules over
additively idempotent semiring having only two trivial strong one-side
ideals (Proposition 4.3) and, based on this characterization, establish some
fundamental facts about e-injective semimodules and their relationship with
injective semimodules over additively idempotent division semirings
(Theorems 4.4 and 4.5) that constitute one of the main goals and results of
the paper. Also, we demonstrate (Proposition-Example 4.6) that the concepts
of `injectivity' and `e-injectivity' for semimodules over semirings, in
general, are different.

In Section 5, we characterize semirings all of whose semimodules are
e-injective (Theorem 5.3), quasi-Frobenius rings in terms of projective and
e-injective semimodules (Theorem 5.4), and one-sided Noetherian rings in
terms of direct sums of e-injective semimodules (Theorem 5.5). These results
are the `e-injective' versions of Theorems 3.4, 3.5, and 3.6 of \cite%
{il'in:dsoisadpopsos}, respectively, and they also are the main results and
another main goal of the paper.

Simple semirings, which are the subject of another important area of
research in the theory of semirings, have quite interesting and promising
applications in various fields (for example, in constructing novel semigroup
actions for a potential use in public-key cryptosystems \cite{mmr:pkcbosa}).
In contrast to the varieties of groups and rings, research on simple
semirings started only recently, and therefore not much on the subject is
known (for some recent results on simple semirings, one may consult \cite%
{bshhurtjankepka:scs}, \cite{monico:ofcss}, \cite{zumbr:cofcsswz}, \cite%
{knt:mosssparp}, \cite{kz:fsais}, \cite{knz:ososacs}, \cite{kkn:srwz}, \cite%
{kn:srwlmae}, and \cite{knz:solpawcics}). While complete characterizations
of commutative and finite simple semirings have been given in \cite%
{bshhurtjankepka:scs} and \cite{kz:fsais}, respectively, the classification
of simple infinite semirings remains an important unresolved problem (see
\cite{knt:mosssparp}, \cite{knz:ososacs}, \cite{kkn:srwz}, and \cite%
{kn:srwlmae} for some recent results in this regard). In light of this and
as a substantial step towards this endeavor, what constitutes one of the
main goals of the paper as well, in Section 6, we give a complete
description of simple semirings all of whose cyclic (finitely generated)
semimodules are e-injective (Theorem 6.10). Moreover, in Section 6, we show
that a bounded distributive lattice all of whose cyclic (finitely generated)
semimodules are e-injective, in fact, is a finite Boolean algebra (Theorem
6.5) and, applying this result, completely characterize subtractive
semirings all of whose cyclic (finitely generated) semimodules are
e-injective (Theorem 6.7). Certainly, all these theorems belong to the main
results of the paper.

Finally, all notions and facts of categorical algebra, used here without any
comments, can be found in \cite{macl:cwm} or \cite{Bor1994a}; for notions
and facts from semiring theory we refer to \cite{golan:sata}.{\ {\ }}

\section{Preliminaries}

\textbf{2.1 }Recall \cite{golan:sata} that a \textit{semiring}\emph{\/} is a
datum $(S,+,\cdot ,0,1)$ such that the following conditions are
satisfied:\smallskip

(1) $(S,+,0)$ is a commutative monoid with identity element $0$;

(2) $(S,\cdot ,1)$ is a monoid with identity element $1$;

(3) Multiplication is distributive over addition from both sides;

(4) $0s=0=s0$ for all $s\in S$.\smallskip\

A semiring that is not a ring we call a \textit{proper semiring}.\smallskip\
A semiring $S$ is a \textit{division semiring} if $(S\setminus \{0\},\cdot
,1)$ is a group; and $S$ is a \emph{semifield} if it is a commutative
division semiring. Two well-known important examples of semifields are the
so-called \textit{Boolean semifield}\emph{\ }$\mathbf{B}=\{0,1\}$ with $%
1+1=1 $, and the tropical semifield $\mathbf{T}:=(\mathbb{R}\cup \{-\infty
\},\max ,+,-\infty ,0\})$.

As usual, a \textit{left\/ }$S$\textit{-semimodule} over the semiring $S$ is
a commutative monoid $(M,+,0_{M})$ together with a scalar multiplication $%
(s,m)\mapsto sm$ from $S\times M$ to $M$ which satisfies the following
identities for all $s,s^{^{\prime }}\in S$ and $m,m^{^{\prime }}\in M$%
:\medskip

(1) $(ss^{^{\prime }})m=s(s^{^{\prime }}m)$;

(2) $s(m+m^{^{\prime }})=sm+sm^{^{\prime }}$;

(3) $(s+s^{^{\prime }})m=sm+s^{^{\prime }}m$;

(4) $1m=m$;

(5) $s0_{M}=0_{M}=0m$.\medskip

\textit{Right semimodules}\emph{\/} over $S$ and homomorphisms between
semimodules and semirings are defined in the standard manner. And, from now
on, let $\mathcal{M}$ be the variety of commutative monoids, and $\mathcal{M}%
_{S}$ and $_{S}\mathcal{M}$ denote the categories of right and left $S$%
-semimodules, respectively, over a semiring $S$.\smallskip\

\textbf{2.2 }An element $\infty \in M$ of an $S$-semimodule $M$ is \textit{%
infinite} if $\infty +m=\infty $ for every $m\in M$; and $K\leq _{\text{ }%
S}M $ means that $K$ is an $S$-subsemimodule of $M$. Also, we will use the
following subsets of the elements of an $S$-semimodule $M$ :%
\begin{equation*}
\begin{tabular}{lll}
$I^{+}(M)$ & $:=$ & $\{m\in M\,|\,m+m=m\};$ \\
$Z(M)$ & $:=$ & $\{z\in M\mid z+m=m\text{ for some }m\in M\};$ \\
$V(M)$ & $:=$ & $\{m\in M\mid m+m^{\prime }=0\text{ for some }m^{\prime }\in
M\};$%
\end{tabular}%
\end{equation*}%
\smallskip

For a semimodule $M$ $\in |_{S}\mathcal{M}|$, it is obvious that $%
I^{+}(M)\cap V(M)=\{0\}$, and $I^{+}(M)\leq _{\text{ }S}Z(M)\leq _{\text{ }%
S}M$.

A left $S$-semimodule $M$ is \textit{zeroic} (\textit{zerosumfree,
additively idempotent}) if $Z(M)=M$ ($V(M)=0,$ $I^{+}(M)=M$). In particular,
a semiring $S$ is \textit{zeroic} (\textit{zerosumfree, additively idempotent%
}) if $_{S}S$ $\in |_{S}\mathcal{M}|$ is a zeroic (zerosumfree, additively
idempotent) semimodule; and we say that $S$ has an infinite element if $%
_{S}S\in |_{S}\mathcal{M}|$ has it. For example, the \textit{Boolean semiring%
} $\mathbf{B}=\{0,1\}$ is a commutative, zeroic, zerosumfree, additively
idempotent semiring in which $\infty =1$.\smallskip

\textbf{2.3} A subsemimodule $A\leq $ $_{S}M$ of a semimodule $M$ is (%
\textit{strongly}) \textit{subtractive}\emph{\ }if ($m+m^{\prime }\in
A\Rightarrow m,m^{\prime }\in A$) $m,m+m^{\prime }\in A\Rightarrow m^{\prime
}\in A$ for all $m,m^{\prime }\in M$. For each subsemimodule $A\leq $ $_{S}M$%
, the subsemimodule $\overline{A}:=\{m\in M\ |\ m+a\in A$ for some $x\in A\}$
is obviously the smallest subtractive subsemimodule of $_{S}M$ containing
the subsemimodule $A$, and therefore, it is called the \textit{subtractive
closure} of $A$; and clearly that $A$ is subtractive iff $A=\overline{A}$. A
left $S$-semimodule $M$ is \textit{subtractive} iff all $S$-subsemimodules
of $M$ are subtractive. In particular, the semiring $S$ is \textit{left
(right)\ subtractive} iff $S$ is subtractive as a left (right) semimodule
over itself. (For important properties of subtractive semirings and
semimodules, the interested reader may consult \cite{knt:ossss} and \cite%
{knt:mosssparp}.)\medskip \emph{\ }

\textbf{2.4} \textit{Congruences}\emph{\ }on an $S$-semimodule $M$ are
defined in the standard manner, and $\mathrm{Cong}(M)$ denotes the set of
all congruences on $M$. This set is non-empty since it always contains at
least two congruences---the \textit{diagonal congruence} $\vartriangle
_{M}:= $ $\{(m,m)$ $|$ $m\in M$ $\}$ and the \textit{universal congruence}%
\emph{\ }$M^{2}:=\{(m,n)$ $|$ $m,n\in M$ $\}$. Any \ subsemimodule $L\leq $ $%
_{S}M$ of an $S$-semimodule $M$ induces a congruence $\equiv _{L}$ on $M$,
known as the \textit{Bourne congruence}, by setting $m\equiv _{L}m^{\prime }$
iff $m+l=m^{\prime }+l^{\prime }$ for some $l,l^{\prime }\in L$; and $M/L$
denotes the factor $S$-semimodule $M/\equiv _{L}$ having the canonical $S$%
-surjection $\pi _{L}:M\longrightarrow M/L$. Following \cite%
{bshhurtjankepka:scs}, a semiring $S$ is \textit{congruence-simple} iff the
only congruences on $S$ are the diagonal $\vartriangle _{S}$ and the
universal $S^{2}$; and $S$ is \textit{ideal-simple} iff $S$ has exactly two
ideals (namely $0$ and $S$). Note that these notions are not the same (see,
\textit{e.g.}, \cite[Examples 3.8]{knz:ososacs}). Moreover, we say that a $S$
is \textit{simple} iff it is simultaneously congruence-simple and
ideal-simple. The classification of simple infinite semirings remains an
important unresolved problem (see \cite{knt:mosssparp}, \cite{knz:ososacs},
\cite{kkn:srwz}, and \cite{kn:srwlmae} for recent related results).\medskip

\textbf{2.5 }As usual (see, for example, \cite[Chapter 17]{golan:sata}), if $%
S$ is a semiring, then in the category $_{S}\mathcal{M}$, a \textit{free}
(left) semimodule $\sum_{i\in I}S_{i},S_{i}\cong $ $_{S}S$, $i\in I$, with a
basis set $I$ is a direct sum (a coproduct) of $|I|$ copies of $_{S}S$; a
semimodule $P\in |_{S}\mathcal{M}|$ is \textit{projective} if it is a
retract of a free semimodule. Following \cite{kat:tpaieosoars} and \cite%
{kat:thcos}, in which there were introduced and considered in detail the
tensor product bifunctors $-\otimes _{S}-$ $:\mathcal{M}_{S}\times $ $_{S}%
\mathcal{M}\longrightarrow \mathcal{M}$, a semimodule $F\in |_{S}\mathcal{M}%
| $ is called (\textit{mono-)flat} iff the functor $-\otimes _{S}F:\mathcal{M%
}_{S}\longrightarrow \mathcal{M}$ preserves (monomorphisms) finite limits;
and the latter is equivalent to that $F$ is a filtered (directed) colimit of
finitely generated free (projective) semimodules. A semimodule $M\in |_{S}%
\mathcal{M}|$ is \textit{finitely generated\ }(\textit{cyclic}) iff $M$ is a
homomorphic image of a free left $S$-semimodule with a finite basis (a
homomorphic image of $_{S}S$); a semimodule $M\in |_{S}\mathcal{M}|$ is
\textit{injective} if for any monomorphism $\mu :A\rightarrowtail B$ of left
$S$-semimodules $A$ and $B\ $and every homomorphism $f\in $ $_{S}\mathcal{M(}%
A,M)$, there exists a homomorphism $\widetilde{f}\in $ $_{S}\mathcal{M(}B,M)$
such that $\widetilde{f}\mu =f$.\smallskip

\textbf{2.6 }For the reader's convenience, we also recall some fundamental
notions and facts about Morita equivalence of semirings (see, \emph{e.g.,}
\cite{kat:thcos} and \cite{kn:meahcos}) that we will use in sequence. Recall
that two semirings $S$ and $T$ are said to be \textit{Morita equivalent} iff
the semimodule categories $_{S}\mathcal{M}$ and $_{T}\mathcal{M}$ are
equivalent, \textit{i.e.}, there exist two functors $F:$ $_{S}\mathcal{M}%
\longrightarrow $ $_{T}\mathcal{M}$ and $G:$ $_{T}\mathcal{M}\longrightarrow
$ $_{S}\mathcal{M}$ and natural isomorphisms $\eta :GF\longrightarrow
Id_{_{S}\mathcal{M}}$ and $\xi :FG\longrightarrow Id_{_{T}\mathcal{M}}$.
Following \cite{kn:meahcos}, a right semimodule $P_{S}$ is said to be a
\textit{generator} for the category $\mathcal{M}_{S}$ of right $S$%
-semimodules iff $S_{S}$ is a retract of a finite direct sum $\oplus _{i}P$
of $P_{S}$; and that $P_{S}$ is said to be a \textit{progenerator} for $%
\mathcal{M}_{S}$ iff $P_{S}$ is a finitely generated projective generator.
Then, by \cite[Theorems 4.5 and 4.12]{kn:meahcos}, $F\simeq P\otimes _{S}-$
for some $(T,S)$-bisemimodule $P$ such that $P_{S}$ is a progenerator, $%
P^{\ast }:=\mathrm{Hom}_{S}(P_{S},S_{S})$ is a progenerator in $\mathcal{M}%
_{T}$, $T\simeq End(P_{S})$ as semirings and $G\simeq P^{\ast }\otimes _{T}-$%
.\medskip

\textbf{2.7 \ }For any left $S$-semimodule $M$, there exists the left $R$-%
\textit{module of differences} $D(M)$ of $M$ \cite[Chapter 16]{golan:sata}
(see also \cite[p. 5083]{kn:orosarp}) defined as the factor semimodule of
the left $S$-semimodule $M\times M$ with respect to the subsemimodule $%
W=\{(m,m)$ $|$ $m\in M$ $\}\subseteq $ $M\times M$, \textit{i.e.}, $%
D(M):=(M\times M)/W$. In fact, the semimodule $D(M)$ is a left $S$-module
since for any $(m,m^{\prime })\in M\times M$ in $D(M)$ one has $\overline{%
(m,m^{\prime })}+\overline{(m^{\prime },m)}=\overline{(0,0)}$. Also, there
exists the canonical $S$-homomorphism $\xi _{M}:$ $M$ $\longrightarrow D(M)$
given by $m\longmapsto \overline{(m,0)}$. In the case when $M$ is a
cancellative semimodule, $\xi _{M}$ is injective, and therefore, we can
consider the elements $\overline{(m,0)}$ and $m$ to be the same and any
element $\overline{(m,m^{\prime })}\in $ $D(M)$ to be the \textquotedblleft
difference\textquotedblright\ of the elements $\overline{(m,0)}$ and $%
\overline{(m^{\prime },0)}$, \textit{i.e.}, $D(M)=\{m-m^{\prime
}\,|\,m,m^{\prime }\in M\}$. In particular, the left $S$-module of
differences $D(S)$ of the regular semimodule $_{S}S$ can be considered as a
ring --- the \textit{ring of differences} of $S$ \cite[Chapter 8, p. 101]%
{golan:sata} --- with the operation of multiplication defined for all $%
a,b,c,d\in S$ by $\overline{(a,b)}\overline{(c,d)}=\overline{(ac+bd,ad+cb)}$%
; and if $S$ is a semiring, then the ring of differences $D(S)$ is also a
semiring with the identity $\overline{(1,0)}$. Moreover, it is easy to see
that $D(M)$ becomes a left $D(S)$-module with $\overline{(a,b)}$ $\overline{%
(m_{1},m_{2})}=\overline{(am_{1}+bm_{2},am_{2}+bm_{1})}$ for all $a,b\in S$
and $m_{1},m_{2}\in M$.\medskip

\textbf{2.8 \ }Any homomorphism $f:M\longrightarrow N$ of left $S$%
-semimodules induces \textit{kernel congruence} $\equiv _{f}$ on $M\ $such
that for any $m,m^{\prime }\in M$, we have $m\equiv _{f}m^{\prime }$ iff $%
f(m)=f(m^{\prime })$, as well as the following subsemimodules:

$\ \ \ \ \ \ \ \ \ \ \ \ \ \ \ \ \ \ \ \ \
\begin{tabular}{lll}
$\mathrm{Ker}(f)$ & $:=$ & $\{m\in M\mid f(m)=0\}$; \\
$im(f)$ & $:=$ & $\{f(m)\mid m\in M\}$; \\
$Im(f)$ & $:=$ & $\{n\in N\text{ }\mid n+f(m)=f(m^{\prime })\text{ for some }%
m,m^{\prime }\in M\}$,%
\end{tabular}%
$

\noindent called the \textit{kernel}, the \textit{image}, and the \textit{%
extended image} of $f$, respectively. Notice that $Im(f)=\overline{im(f)}$
and that $\equiv $ $_{\mathrm{Ker}(f)}$ $\subseteq $ $\equiv _{f}$.
Moreover, one can easily see that $N/im(f)=N/Im(f)$ and $M/\equiv _{f}$ $%
\simeq im(f)$. However, in general, in our non-abelian setting, as the
following example shows, the $S$-semimodules $M/\mathrm{Ker}(f)$ and $im(f)$
are not necessarily isomorphic: Indeed, the tropical semifield $\mathbf{T}:=(%
\mathbb{R}\cup \{-\infty \},\max ,+,-\infty ,0\})$ contains the Boolean
semifield $\mathbf{B}$ as a subsemiring, whence $\mathbf{T}$ can be
considered as a $\mathbf{B}$-semimodule in a canonical way, and let $f:%
\mathbf{T}\longrightarrow \mathbf{B}$ be the $\mathbf{B}$-homomorphism map
given by%
\begin{equation*}
f(x)=\left\{
\begin{array}{ccc}
0, &  & x=-\infty \\
&  &  \\
1, &  & x\neq -\infty%
\end{array}%
\right. \text{;}
\end{equation*}%
Then, clearly, $f$ is surjective and $\mathrm{Ker}(f)=0$, but $\mathbf{T}/%
\mathrm{Ker}(f)\simeq \mathbf{T}\ncong \mathbf{B}=im(f)$.

\section{e-Injectivity (-projectivity) of semimodules}

For the readers' convenience, briefly remind some general notions for
categories with zero morphisms\ (see, for example, \cite[Chapters 7, 8, 13]%
{schubert:c}) in the context of semimodule categories $_{S}\mathcal{M}$. A
\textit{kernel} $(K,k)$ of a homomorphism $f:A\longrightarrow B$ of left $S$%
-semimodules is a homomorphism $k:K\longrightarrow A$ such that \ (i) $fk=0$%
, (ii) for every homomorphism $x:X\longrightarrow A$ with $fx=0$ there
exists exactly one homomorphism $i:X\longrightarrow K$ such that $x=ki$; and
in this case, we write $(K,k)=\ker f$ . By dualizing, one comes up to the
concept of a \textit{cokernel} $(k,K)$, $k:B\longrightarrow K$ of $%
f:A\longrightarrow B$, and\ $(k,K)=co\ker f$. Then, it is obvious (or it is
just readily follows from (co)completeness of semimodule categories\ $_{S}%
\mathcal{M}$ \cite[Sections 5.1, 5.2]{macl:cwm}) that there exist $\ker f$ \
and $co\ker f$ for every homomorphism $f\in $ $_{S}\mathcal{M}(A,B):=\mathrm{%
Hom}_{S}(A,B)$. Furthermore, following, for example, \cite[Sect. II.6.2]%
{gelman:moha} or \cite[Sect.13.2]{schubert:c} and without loss of
generality, we can define \textit{a short exact sequence} in $_{S}\mathcal{M}
$ as a sequence%
\begin{equation*}
0\longrightarrow A\overset{f}{\longrightarrow }B\overset{g}{\longrightarrow }%
C\longrightarrow 0\text{ \ \ \ \ \ \ \ \ \ \ \ \ \ }(\ast )
\end{equation*}%
of semimodules $A,B,C\in |_{S}\mathcal{M}|$ and homomorphisms $f$ \ and $g$
such that $(A,f)=\ker g$ and $(g,C)=co\ker f$.

Let $F:$ $_{S}\mathcal{M\longrightarrow }$ $_{R}\mathcal{M}$ ($G:$ $_{S}%
\mathcal{M\longrightarrow }$ $_{R}\mathcal{M}$) be a covariant
(contravariant) functor between the semimodule categories $_{S}\mathcal{M}$
and $_{R}\mathcal{M}$. Then, the functor $F$ ($G$) we say is an \textit{%
exact functor }(or just an \textit{e-functor}) if for any exact sequence (*)
in $_{S}\mathcal{M}$ the sequence
\begin{equation*}
0\longrightarrow F(A)\overset{F(f)}{\longrightarrow }F(B)\overset{F(g)}{%
\longrightarrow }F(C)\longrightarrow 0\text{ (}0\longleftarrow G(A)\overset{%
G(f)}{\longleftarrow }G(B)\overset{G(g)}{\longleftarrow }G(C)\longleftarrow 0%
\text{)}
\end{equation*}%
is exact in $_{R}\mathcal{M}$ as well. The following observations will prove
to be useful.\medskip

\noindent \textbf{Proposition 3.1 }\textit{(1) A direct sum} $\oplus
_{i}F_{i}$ \textit{of covariant functors }$F_{i}:$ $_{S}\mathcal{%
M\longrightarrow }$ $_{R}\mathcal{M}$\textit{,} $i\in I$\textit{, is an
e-functor iff each summand }$F_{i}$\textit{\ is an e-functor;}

\textit{(2)} \textit{A direct product} $\Pi _{i\in I}$ $G_{i}$ \textit{of
contravariant functors} $G_{i}:$ $_{S}\mathcal{M\longrightarrow }$ $_{R}%
\mathcal{M}$\textit{, }$i\in I$\textit{, is an e-functor iff each factor }$%
G_{i}$\textit{\ is an e-functor;}

\textit{(3) A retract of an e-functor is an e-functor as well.\medskip }

\noindent \textit{Proof} \ (1). It immediately follows from the observation
that for an exact sequence $(\ast )$, the sequence%
\begin{equation*}
0\longrightarrow \oplus _{i}F_{i}(A)\overset{\oplus _{i}F_{i}(f)}{%
\longrightarrow }\oplus _{i}F_{i}(B)\overset{\oplus _{i}F_{i}(g)}{%
\longrightarrow }\oplus _{i}F_{i}(C)\longrightarrow 0\text{ }
\end{equation*}%
is, in fact, a `direct sum'%
\begin{equation*}
0\longrightarrow \oplus _{i}(F_{i}(A))\overset{\oplus _{i}(F_{i}(f))}{%
\longrightarrow }\oplus _{i}(F_{i}(B))\overset{\oplus _{i}(F_{i}(g))}{%
\longrightarrow }\oplus _{i}(F_{i}(C))\longrightarrow 0
\end{equation*}%
of the exact sequences%
\begin{equation*}
0\longrightarrow F_{i}(A)\overset{F_{i}(f)}{\longrightarrow }F_{i}(B)\overset%
{F_{i}(g)}{\longrightarrow }F_{i}(C)\longrightarrow 0\text{, \ \ }i\in I%
\text{,}
\end{equation*}%
and therefore, is an exact sequence itself.

(2). It is established in the same fashion as (1).

(3). Clearly, it is enough to consider only a \textquotedblleft
covariant\textquotedblright\ case. So, let $Q$ be a retract of an exact
covariant functor $F:$ $_{S}\mathcal{M\longrightarrow }$ $_{R}\mathcal{M}$,
\textit{i.e.}, there exist natural transformations of functors $\mu
:Q\longrightarrow F$ and $\pi :F\longrightarrow Q$ such that $\pi \mu =1_{Q}$%
. Then, a \textquotedblleft diagram chase\textquotedblright\ of a
commutative diagram for an exact sequence $(\ast )$%
\begin{equation*}
\begin{tabular}{lllllllll}
$0$ & $\longrightarrow $ & $Q(A)$ & $\overset{Q(f)}{\longrightarrow }$ & $%
Q(B)$ & $\overset{Q(g)}{\longrightarrow }$ & $Q(C)$ & $\longrightarrow $ & $%
0 $ \\
$\downarrow $ &  & $\downarrow \mu _{A}$ &  & $\downarrow \mu _{B}$ &  & $%
\downarrow \mu _{C}$ &  & $\downarrow $ \\
$0$ & $\longrightarrow $ & $F(A)$ & $\overset{F(f)}{\longrightarrow }$ & $%
F(B)$ & $\overset{F(g)}{\longrightarrow }$ & $F(C)$ & $\longrightarrow $ & $%
0 $ \\
$\downarrow $ &  & $\downarrow \pi _{A}$ &  & $\downarrow \pi _{B}$ &  & $%
\downarrow \pi _{C}$ &  & $\downarrow $ \\
$0$ & $\longrightarrow $ & $Q(A)$ & $\overset{Q(f)}{\longrightarrow }$ & $%
Q(B)$ & $\overset{Q(g)}{\longrightarrow }$ & $Q(C)$ & $\longrightarrow $ & $%
0 $%
\end{tabular}%
\end{equation*}

\noindent leads us to the exactness of the sequence $0\longrightarrow Q(A)%
\overset{Q(f)}{\longrightarrow }Q(B)\overset{Q(g)}{\longrightarrow }%
Q(C)\longrightarrow 0$. \textit{\ \ \ \ \ \ }$_{\square }\medskip $

In the case when a semiring $S$ is a ring, the category $_{S}\mathcal{M}$
becomes an abelian category of modules over a ring $S$. Then, there are two,
at the first glance, different approaches to the concepts of injectivity and
projectivity of modules: the first one is an `universal algebra' approach
presented in 2.5, and the second --- a `homological algebra' approach that
in more general setting of semimodules over semirings is given by the
following definition.\medskip

\noindent \textbf{Definition 3.2 }A semimodule $M\in |_{S}\mathcal{M}|$ over
a semiring $S$ is \textit{exactly}-\textit{injective }(\textit{%
exactly-projective}), or shortly \textit{e}-\textit{injective }(\textit{%
e-projective}), if $_{S}\mathcal{M}(-,M):=$

\noindent $\mathrm{Hom}_{S}(-,M):$ $_{S}\mathcal{M\longrightarrow }$ $%
\mathcal{M}$ ($_{S}\mathcal{M}(M,-):=\mathrm{Hom}_{S}(M,-):$ $_{S}\mathcal{%
M\longrightarrow }$ $\mathcal{M}$) is an exact contravariant (covariant)
functor.\medskip

Perhaps one of the most important achievements of homological algebra are
due to the fact that for abelian categories of modules $_{S}\mathcal{M}$
over a ring $S$\ the both --- universal algebra and homological algebra ---
approaches lead to the identical classes of modules (see, for example, \cite[%
Sections 1.2A and 1.3A]{lam:lomar} or \cite[Sect. II.6.9]{gelman:moha}).
However, as we will see later on, in generally non-abelian setting of
categories of semimodules $_{S}\mathcal{M}$ over a semiring $S$, these two
approaches lead to two different classes of semimodules. As was mentioned
earlier, in non-abelian (even in non-additive) settings, there have been
obtained a good number of quite interesting and important results connected
with the concepts of `projectivity' and `injectivity' based on the universal
algebra approach. In contrast to this, in the present paper, we initiate
investigations related to injective and projective semimodules defined by
using the second approach, \textit{i.e.}, from the homological algebra point
of view. Thus, from Proposition 3.1 we right away obtain\medskip

\noindent \textbf{Corollary 3.3 }\textit{(1) A direct sum} $\oplus _{i}M_{i}$
\textit{of e-projective semimodules }$M_{i}$ $\in $ $|_{S}\mathcal{M}|$%
\textit{,} $i\in I$\textit{, is an e-projective semimodule iff each summand }%
$M_{i}$\textit{\ is an e-projective semimodule;}

\textit{(2)} \textit{A direct product} $\Pi _{i\in I}$ $M_{i}$ \textit{of
e-injective semimodules }$M_{i}$ $\in $ $|_{S}\mathcal{M}|$\textit{,} $i\in
I $\textit{, is an e-injective semimodule iff each factor }$M_{i}$\ \textit{%
is an} \textit{e}-\textit{injective} \textit{semimodule;}

\textit{(3) A retract of an e-injective (e-projective) semimodule is an
e-injective (e-projective) semimodule as well.\medskip }

Now, let $_{S}\mathcal{M}_{T}$ be a category of bisemimodules over semirings
$T$ and $S$ (see, \textit{e.g.}, \cite[Sect. 3]{kat:thcos}) and $P\in |_{S}%
\mathcal{M}_{T}|$. Then for any semimodule $M$ $\in $ $|_{S}\mathcal{M}|$,
on the set of homomorphisms $\mathrm{Hom}_{S}(_{S}P,_{S}M):=$ $_{S}\mathcal{M%
}(P,M)$ there exists a natural structure of a left $T$-semimodule defined as
follows: $(tf)(p)=f(pt)$ for any $f\in \mathrm{Hom}_{S}(_{S}P,_{S}M)$, $p\in
P$, and $t\in T$. Moreover, there is a functor $\mathrm{Hom}_{S}(_{S}P,-):$ $%
_{S}\mathcal{M\longrightarrow }$ $_{T}\mathcal{M}$. The next result provides
us with main tools\ for constructing examples of (e-)injective semimodules
by using known (e-)injective semimodules over one semiring to produce
(e-)injective semimodules over another. It includes and/or generalizes the
well-known \textquotedblleft Injective Producing Lemma\textquotedblright\
\cite[Lemma 1.3.5]{lam:lomar} and \cite[Proposition 17.25]{golan:sata}%
.\medskip

\noindent \textbf{Proposition 3.4 }\textit{(1) Let for a semimodule} $P\in
|_{S}\mathcal{M}_{T}|$\textit{\ the semimodule} $P_{T}\in |\mathcal{M}_{T}|$
\textit{be mono-flat and a semimodule }$M$ $\in $ $|_{S}\mathcal{M}|$\textit{%
\ an injective left }$S$\textit{-semimodule, then the semimodule }$\mathrm{%
Hom}_{S}(_{S}P,_{S}M)\in |_{T}\mathcal{M}|$ \textit{is an injective left }$T$%
\textit{-semimodule as well;}

\textit{(2) Let for a semimodule }$P\in |_{S}\mathcal{M}_{T}|$\textit{\ the
functor} $P_{T}$ $\otimes -:$ $_{T}\mathcal{M\longrightarrow }_{S}\mathcal{M}
$ \textit{preserve short exact sequences and a semimodule }$M$ $\in $ $|_{S}%
\mathcal{M}|$ \textit{an e-injective left }$S$\textit{-semimodule, then the
semimodule} $\mathrm{Hom}_{S}(_{S}P,_{S}M)\in |_{T}\mathcal{M}|$ \textit{is
an e-injective left }$T$\textit{-semimodule as well.\medskip }

\noindent \textit{Proof} \ (1). It is easy to see that a semimodule $%
_{T}I\in |_{T}\mathcal{M}|$ is injective iff the contravariant functor $%
\mathrm{Hom}_{T}(-,_{T}I):$ $_{T}\mathcal{M\longrightarrow M}$ moves any
monomorphism $\mu :A\rightarrowtail B$ of left $T$-semimodules $A$ and $B$
to the surjection $\mathrm{Hom}_{T}(\mu ,_{T}I):$ $\mathrm{Hom}_{T}(B,_{T}I)$
$\mathcal{\longrightarrow }$ $\mathrm{Hom}_{T}(A,_{T}I)$ in $\mathcal{M}$.
Therefore, we need only to show that for any injective semimodule $M$ $\in $
$|_{S}\mathcal{M}|$ and any monomorphism $\mu :A\rightarrowtail B$ in $_{T}%
\mathcal{M}$,
\begin{eqnarray*}
\mathrm{Hom}_{T}(\mu ,\mathrm{Hom}_{S}(_{S}P,_{S}M)) &:&\mathrm{Hom}_{T}(B,%
\mathrm{Hom}_{S}(_{S}P,_{S}M))\text{ }\mathcal{\longrightarrow }\text{ } \\
&&\mathrm{Hom}_{T}(A,\mathrm{Hom}_{S}(_{S}P,_{S}M))
\end{eqnarray*}%
is a surjection in $\mathcal{M}$. However, by \cite[Theorem 3.3]{kat:thcos}
(or \cite[Theorem 3.6]{kat:tpaieosoars}), the functor $P_{T}$ $\otimes -:$ $%
_{T}\mathcal{M\longrightarrow }_{S}\mathcal{M}$ is a left adjoint to the
functor $\mathrm{Hom}_{S}(_{S}P,-):$ $_{S}\mathcal{M\longrightarrow }$ $_{T}%
\mathcal{M}$. Whence, from the commutative diagram%
\begin{equation*}
\begin{tabular}{lll}
$\mathrm{Hom}_{S}(_{S}P\otimes A,_{S}M)$ & $\simeq $ & $\mathrm{Hom}_{T}(A,%
\mathrm{Hom}_{S}(_{S}P,_{S}M))$ \\
$\downarrow \mathrm{Hom}_{S}(1_{_{S}P}\otimes \mu ,1_{_{S}M})$ &  & $%
\uparrow \mathrm{Hom}_{T}(\mu ,\mathrm{Hom}_{S}(_{S}P,_{S}M))$ \\
$\mathrm{Hom}_{S}(_{S}P\otimes B,_{S}M)$ & $\simeq $ & $\mathrm{Hom}_{T}(B,%
\mathrm{Hom}_{S}(_{S}P,_{S}M))$%
\end{tabular}%
\end{equation*}%
and the mono-flatness of the semimodule $P_{T}$, we have that $\mathrm{Hom}%
_{S}(_{S}P,_{S}M)\in |_{T}\mathcal{M}|$ is an injective semimodule.

(2). In the same fashion as it was done in (1) and assuming that $P_{T}$ $%
\otimes -:$ $_{T}\mathcal{M\longrightarrow }_{S}\mathcal{M}$ preserves short
exact sequences, one gets the e-injectivity of

\noindent $\mathrm{Hom}_{S}(_{S}P,_{S}M)\in |_{T}\mathcal{M}|$. \textit{\ \
\ \ \ \ }$_{\square }\medskip $

\noindent\ \textbf{Corollary 3.5 }(cf. \cite[Corollary 1.3.6B]{lam:lomar}
and \cite[Proposition 17.25]{golan:sata}) \textit{Let }$f:S\longrightarrow T$%
\textit{\ be a semiring homomorphism. Then the functor }$\mathrm{Hom}%
_{S}(_{S}T,-):$ $_{S}\mathcal{M\longrightarrow }$ $_{T}\mathcal{M}$ \textit{%
preserves (e-)injective semimodules.\medskip }

\noindent \textit{Proof} \ Taking into consideration that by \cite[%
Proposition 3.8]{kat:tpaieosoars} there is the functor isomorphism $T_{T}$ $%
\otimes -\simeq Id:$ $_{T}\mathcal{M\longrightarrow }\ \mathcal{M}$, and
therefore, $T_{T}$ $\otimes -$ is mono-flat and preserves all (co)limits as
well, the statements follow right away from Proposition 3.4. \textit{\ \ \ \
\ \ }$_{\square }\medskip $

We conclude this section with the following two remarks.\medskip

\noindent \textbf{Remark 3.6 }By dualization, the interested reader may
easily obtain the corresponding analogs of the statements of Proposition 3.4
and Corollary 3.5 regarding (e-)projective semimodules.\medskip

\noindent \textbf{Remark 3.7 }It should be mentioned that the concept of
e-injectivity for semimodules, in some way and using different terminology,
have been earlier considered in \cite{Tak82}, \cite{Tak:eos}, \cite%
{Abu2014-a}, and \cite{Abu2014-b}, where the authors heavily used the
obvious fact that in categories of semimodules a monomorphism $\mu
:A\rightarrowtail B$ is a kernel of its cokernel iff the subsemimodule $%
A\leq $ $B$ is a subtractive one. Moreover, it is easy to see that a
semimodule $_{T}I\in |_{T}\mathcal{M}|$ is e-injective iff the contravariant
functor $\mathrm{Hom}_{T}(-,_{T}I):$ $_{T}\mathcal{M\longrightarrow M}$
moves any monomorphism $\mu :A\rightarrowtail B$ with $\mu (A)$ being a
subtractive subsemimodule of a left $T$-semimodule $B$ to the surjection $%
\mathrm{Hom}_{T}(\mu ,_{T}I):$ $\mathrm{Hom}_{T}(B,_{T}I)$ $\mathcal{%
\longrightarrow }$ $\mathrm{Hom}_{T}(A,_{T}I)$ in $\mathcal{M}$ and $(%
\mathrm{Hom}_{T}(\mu ,_{T}I)$, $\mathrm{Hom}_{T}(A,_{T}I))=co\ker (\mathrm{%
Hom}_{T}(\pi ,_{T}I))$, where $\pi:\, B\longrightarrow B/\mu(A)$ is the natural surjection;
in other words, $_{T}I\in |_{T}\mathcal{M}|$ is e-injective iff, for any
$T$-homomorphism $\varphi:\, A\longrightarrow I$, there exists a $T$-homomorphism
$\psi:\, B\longrightarrow I$ such that $\varphi=\psi\mu$ and, for any $T$-homomorphisms
$\psi_1$, $\psi_2:$ $B$ $\mathcal{\longrightarrow }$ $I$, the equality
$\psi_1\mu=\psi_2\mu$ implies $\psi_1+\chi_1=\psi_2+\chi_2$ for suitable $\chi_1$
and $\chi_2$ for which $\chi_1\mu=0=\chi_2\mu$. And for that
reason, the terminology used in those works was certainly influenced by
these facts. Our approach and \textquotedblleft ideology\textquotedblright\
are slightly different.

\section{ e-Injective semimodules over additively idempotent semirings}

In this section, we consider the structure of e-injective semimodules as
well as relations between e-injective and injective semimodules over
additively idempotent semirings. To do so, we need first to establish some
useful facts.\medskip

\noindent \textbf{Proposition 4.1 }\textit{Let }$S$\textit{\ be a
zerosumfree semiring and} $M$ $\in $ $|_{S}\mathcal{M}|$ \textit{an
e-injective }$S$\textit{-semimodule. Then there exists an element }$z\in M$%
\textit{\ such that }$m+z=z$\textit{\ for all }$m\in M$\textit{, and in
addition, }$sz=z$\textit{\ for every }$s\in S\backslash \{0\}$\textit{\ if }$%
S$\textit{\ is an entire semiring.\medskip }

\noindent \textit{Proof} \ For an e-injective semimodule $M\in $ $|_{S}%
\mathcal{M}|$ consider the relation $\sim $ on the left $S$-semimodule $%
S\times M\in $ $|_{S}\mathcal{M}|$ defined as follows: $(s_{1},m_{1})\sim
(s_{2},m_{2})$ iff $s_{1}=s_{2}$ and there are $m_{i}^{\prime
},m_{i}^{\prime \prime }\in M$, $x_{i},x_{i}^{\prime }\in S$ $(i=1,2,...,n)$
such that $x_{i}+x_{i}^{\prime }=s_{1}$ and $m_{1}+\Sigma
_{i=1}^{n}x_{i}m_{i}^{\prime }=m_{2}+\Sigma _{i=1}^{n}x_{i}m_{i}^{\prime
\prime }$ for some natural number $n$. As was shown in \cite[Propositions
1.6 and 1.7]{il'in:dsoisadpopsos}, $\sim $ is an $S$-congruence on $S\times
M\in $ $|_{S}\mathcal{M}|$ and there exists an $S$-monomorphism $\alpha
:M\longrightarrow S\times M/{\sim ,}$ $m\mapsto \lbrack 0_{S},m]$ such that $%
[0_{S},m]+[1,0_{M}]=[1,0_{M}]$ for all $m\in M$, respectively. Obviously, $%
\alpha (M)$ is a subtractive subsemimodule of the semimodule $S\times M/{%
\sim }$. Therefore as well as taking into considertation Remark 3.7, there
exists an $S$-homomorphism $\beta :S\times M/{\sim }\longrightarrow M$ such
that $\beta\alpha =id_{M}.$ Let $z:=\beta ([1,0_{M}])$; and as \ it
was shown in \cite[Proposition 1.7]{il'in:dsoisadpopsos}, $m+z=z$ for all $%
m\in M$.

If in addition, $S$ is an entire semiring, it is easy to see that for the
natural embedding $\mu :M\rightarrowtail Ext(M)$, where $Ext(M):=M\cup
\{\infty \}$ with $0\infty =0$ and $s\infty =\infty $ for all $s\in
S\backslash \{0\}$ and $m+\infty =\infty $ for all $m\in Ext(M)$, there is a
subtractive semimodule $\mu (M)\leq $ $Ext(M)$. Since $_{S}M$ is
e-injective, by Remark 3.7, there exists $\gamma :Ext(M)\longrightarrow M$
such that $\gamma \mu =id_{M}$, and clearly $z=\gamma (\infty )$ and $%
sz=s\gamma (\infty )=\gamma (s\infty )=\gamma (\infty )=z$ for every $s\in
S\backslash \{0\}$. \textit{\ \ \ \ \ \ }$_{\square }\medskip $

\noindent \textbf{Lemma 4.2 }\textit{If a zerosumfree semiring }$S$\textit{\
has only two trivial strongly subtractive left (right) ideals, then }$S$%
\textit{\ is entire.\medskip }

\noindent \textit{Proof} \ Indeed, suppose that $ab=0$ for some $a,b\in S$
and $b\neq 0.$ Since $S$ is zerosumfree, the left ideal%
\begin{equation*}
(0:\text{ }_{S}b):=\{s\in S\,|\,sb=0\}
\end{equation*}%
is strongly subtractive and $1\notin (0:$ $_{S}b)$. Hence, $(0:$ $_{S}b)=0$
and $a=0$. \textit{\ \ \ \ \ \ }$_{\square }\medskip $

The next observation, complementing Proposition 4.1, completely describes
e-injective semimodules over additively idempotent semirings with only two
trivial one-sided strongly subtractive ideals.\medskip

\noindent \textbf{Proposition 4.3 }\textit{A left }$S$\textit{-semimodule} $%
M $ $\in $ $|_{S}\mathcal{M}|$ \textit{over an additively idempotent
semiring }$S$\textit{\ with only two trivial left strongly subtractive
ideals is e-injective iff it possesses the infinite element }$\infty \in $%
\textit{\ }$M$\textit{\ such that }$s\infty =\infty $\textit{\ for all }$%
s\in S\backslash \{0\}$\textit{.}$\medskip $

\noindent \textit{Proof} \ $\Longrightarrow $. By Lemma 4.2, $S$ is an
additively idempotent, zerosumfree, and entire, semiring. Whence, the result
follows from Proposition 4.1.

$\Longleftarrow $. Let $\mu :A\rightarrowtail B$ be the inclusion map from a
subtractive subsemimodule $_{S}A$ to a semimodule $_{S}B$ in $_{S}\mathcal{M}$.
For $S$ is additively idempotent, the additive reducts of all $S$%
-semimodules are idempotent monoids as well. Moreover, $_{S}A$ is a strongly
subtractive subsemimodule of $_{S}B$: Indeed, if $x+y\in A$ for some $x,y\in
B$, then $x+(x+y)=x+y\in A$ and $y+(x+y)=x+y\in A$, and hence, $x,y\in A$.

Next, for an $S$-semimodule homomorphism $f:A\longrightarrow M$, define a map
$g:B\longrightarrow M$ as follows:%
\begin{equation*}
g(b)=\left\{
\begin{array}{lcl}
f(b), &  & b\in A \\
&  &  \\
\infty , &  & b\notin A%
\end{array}%
\right. \text{.}
\end{equation*}%
It is obvious that $f=g\mu $, and we claim that actually $g$ is an $S$%
-semimodule homomorphism, too.

Indeed, for any $b,c\in B$, if $b+c\in A$, then $b,c\in A$ since $A$ is a
strongly subtractive subsemimodule of $B$, and therefore, $%
g(b+c)=f(b+c)=f(b)+f(c)=g(b)+g(c)$; if $b+c\notin A,$ then we have the
following two cases to consider:

\emph{Case 1: }$b\in A$ and $c\notin A$. Then, $g(b+c)=\infty =f(b)+\infty
=g(b)+g(c)$;

\emph{Case 2:} $b\notin A$ and $c\notin A$. Then, $g(b+c)=\infty =\infty
+\infty =g(b)+g(c)$.

Now let $b\in B$ and $s\in S$. For $A\leq _{S}B$ is a strongly subtractive
subsemimodule, the left ideal%
\begin{equation*}
I_{b}:=\{s\in S\,|\,sb\in A\}
\end{equation*}%
is a strongly subtractive left ideal in $S$. Hence, $I_{b}=0$ or $I_{b}=S$;
and $sb\in A$ iff $s=0$ or $b\in A$, or equivalently, $sb\notin A$ iff $%
s\neq 0$ and $b\notin A$. Furthermore, we always have $s\infty =\infty $ for
all $s\in S\backslash \{0\}$. From these observations, we immediately get $%
g(sb)=sg(b)$, and therefore, $g$ is an $S$-semimodule homorphism.

Finally, let $\pi:\, B\longrightarrow B/A$ be the natural surjection. It is
easy to see that there is a homomorphism $\gamma \in \mathrm{Hom}_{S}(B/A,M)$ such
that the image $\gamma ([b])$ of any non-zero element $[b]\in B/A$ is the
infinite element $\infty \in M$. Moreover, it is quite
clear that for any $\alpha,\alpha^{\prime }\in \mathrm{Hom}_{S}(B,M)$ such that
$\alpha\mu=\alpha^{\prime}\mu$, we have \mbox{$\alpha+\gamma\pi$} $=$ \mbox{$\alpha^{\prime}+\gamma\pi$} with
$\gamma\pi\mu=0$. Hence, by Remark~3.7, we get 
$(\mathrm{Hom}_S(\mu,M),\mathrm{Hom}_S(A,M))=co\ker(\mathrm{Hom}_S(\pi,M))$, and $_{S}M$ is an
e-injective $S$-semimodule. \textit{\ \ \ \ \ \ }$_{\square }\medskip $

It is obvious that every division semiring has only two trivial strongly
left/right subtractive ideals. So, as a consequence of Proposition 4.3,
among other, in our view, interesting and important observations, in the
next theorem we obtain a complete characterization of e-injective
semimodules over additively idempotent division semirings.\medskip

\noindent \textbf{Theorem 4.4 }\textit{For an additively idempotent division
semiring }$D$\textit{\ the following statements are true:\medskip }

\textit{(1) A left }$D$\textit{-semimodule }$M$\textit{\ is e-injective iff }%
$M$\textit{\ has an infinite element;}

\textit{(2) Every left }$D$\textit{-semimodule can be embedded in an
e-injective }$D$\textit{-semimodule;}

\textit{(3) Every injective left }$D$\textit{-semimodule is e-injective;}

\textit{(4) Every finite left }$D$\textit{-semimodule is e-injective;}

\textit{(5) The regular }$D$\textit{-semimodule }$_{D}D$ (\textit{and
therefore,} \textit{every finitely generated left }$D$\textit{-semimodule)
is e-injective iff }$D\simeq \mathbf{B}$\textit{.\medskip }

\noindent \textit{Proof} \ So, let $D$ be an additively idempotent, and
therefore, zerosumfree, division semiring.

(1) $\Longrightarrow $. This follows immediately from Proposition~4.3.

$\Longleftarrow $. Let $M$ be a left $D$-semimodule with the infinite
element $\infty $. For every $s\in D\backslash \{0\}$, we have $s\infty
+\infty =\infty $, and hence, \noindent $\infty =\infty +s^{-1}\infty
=s^{-1}\infty ,$ whence $s\infty =\infty $. From this observation and
Proposition~4.3, we conclude that $_{D}M$ is e-injective.

(2) For every left $D$-semimodule $M$, the left $D$-semimodule $%
Ext(M):=M\cup \{\infty \}$ has an infinite element, and therefore by (1), it
is e-injective.

(3) Let $M$ be an injective left $D$-semimodule. By \cite[Proposition 1.7]%
{il'in:dsoisadpopsos}, $M$ has an infinite element $\infty $, whence $M$ is
e-injective by (1).

(4) Let $M=\{m_{1},\cdots ,m_{n}\}$. Then, $\infty :=\Sigma _{i=1}^{n}m_{i}$
is an infinite element of $M\ $and by (1) $_{D}M$ is e-injective.

(5) $\Longrightarrow $. Assume that $_{D}D$ is e-injective. Since $D$ is
zerosumfree and by Proposition 4.1, it has an infinite element $\infty $. It
follows that $\infty =\infty +\infty ^{2}=(\infty +1)\infty =\infty ^{2}$,
and hence, $1=\infty $. Then, for every $d\in D\backslash \{0\}$, we have $%
d^{-1}+1=1$ and $1=d+1=d(1+d^{-1})=d\cdot \ 1=d$. Thus, $D\simeq \mathbf{B}$.

$\Longleftarrow $. This is clear from (1) and (4). \textit{\ \ \ \ \ \ }$%
_{\square }\medskip $

Our next observation extends statements (2) and (3) of Theorem 4.4 to
arbitrary idempotent semirings and actually establishes the e-injective
completeness of categories of $S$-semimodules $_{S}\mathcal{M}$ over
additively idempotent semirings $S$ (the injective completeness of those
categories has been established earlier in \cite[Theorem 4.2]%
{kat:tpaieosoars}), namely:\medskip

\noindent \textbf{Theorem 4.5 }\textit{For an additively idempotent semiring
}$S$\textit{\ the following statements are true:\medskip }

\textit{(1) Every left }$S$\textit{-semimodule can be embedded in an
e-injective }$S$\textit{-semimodule;}

\textit{(2) Every injective left }$S$\textit{-semimodule is
e-injective.\medskip }

\noindent \textit{Proof} \ (1) Since $S$ is additively idempotent, there is
the obvious semiring embedding $\mu :\mathbf{B}\rightarrowtail S$ of the
Boolean semifield $\mathbf{B}$ into a subsemiring $S$. Therefore, any
semimodule $M\in $ $|_{S}\mathcal{M}|$ has the canonical $\mathbf{B}$%
-semimodule structure, \textit{i.e.}, $M\in $ $|_{\mathbf{B}}\mathcal{M}|$.
By (2) of Theorem 4.4, there exists an an e-injective $\mathbf{B}$%
-semimodule $Q$ such that $M\leq $ $_{\mathbf{B}}Q$. Then, taking into
consideration obvious embeddings%
\begin{equation*}
_{S}M\simeq \mathrm{Hom}_{S}(S,M)\leq \mathrm{Hom}_{\mathbf{B}}(S,M)\leq
\mathrm{Hom}_{\mathbf{B}}(S,Q)\text{,}
\end{equation*}%
and applying Corollary 3.5, we have that $\mathrm{Hom}_{\mathbf{B}}(_{S}S,Q)$
is an e-injective left $S$-semimodule.

(2) Let $M\in $ $|_{S}\mathcal{M}|$ be an injective left $S$-semimodule. By
(1), in the category $_{S}\mathcal{M}$ there exists an embedding $\mu
:M\rightarrowtail Q$, where $Q$ is an e-injective $S$-semimodule. Whence, as
$_{S}M$ is an injective $S$-semimodule, there exists an $S$-homomorphism $%
\pi :Q\longrightarrow M$ such that $\pi \mu =id_{M}$, \textit{i.e.}, $_{S}M$
is a retract of $_{S}Q$; and by (3) of Corollary 3.3, we end the proof.
\textit{\ \ \ \ \ \ }$_{\square }\medskip $

Theorems 4.4 and 4.5 show that in the category $_{S}\mathcal{M}$ over an
additively idempotent semiring $S$, the\ class of injective semimodules, $%
Inj(_{S}\mathcal{M)}$, is a subclass of the\ class of e-injective
semimodules, $e$-$Inj(_{S}\mathcal{M)}$. In fact, our next observation shows
that for an additively idempotent division semiring $S$, the subclass $%
Inj(_{S}\mathcal{M)}$ is always a proper one of the class $e$-$Inj(_{S}%
\mathcal{M)}$, namely:\medskip

\noindent \textbf{Proposition-Example 4.6 }\textit{Let }$D$\textit{\ be an
additively idempotent division semiring. Then} $Inj(_{D}\mathcal{M)}$ $%
\subset $ $e$-$Inj(_{D}\mathcal{M)}$\textit{.\medskip }

\noindent \textit{Proof} \ For $D$ is a zerosumfree division semiring, there
exists a surjective semiring homomorphism $\pi :D\longrightarrow \mathbf{B}$%
. By \cite[Lemma 5.2]{kn:meahcos}, the \emph{restriction functor} $\pi
^{\#}: $ $_{\mathbf{B}}\mathcal{M}\longrightarrow $ $_{D}\mathcal{M}$
preserves non-injective semimodules. Let $M\in $ $|_{\mathbf{B}}\mathcal{M}|$
be a non-distributive finite lattice, then $_{\mathbf{B}}M$ is non-injective
$\mathbf{B}$-semimodule \cite[Theorem 4]{fofan:ipoba}, and therefore, $\pi
^{\#}(M)\in $ $|_{D}\mathcal{M}|$ is non-injective $D$-semimodule. However,
by (4) of Theorem 4.4, $\pi ^{\#}(M)\in e$-$Inj(_{D}\mathcal{M)}$. \textit{\
\ \ \ \ \ }$_{\square }\medskip $

In light of the observations above, we finish this section by posting the
following, in our view interesting and perspective, conjecture and
problems.\medskip

\noindent \textbf{Conjecture 1} $Inj(_{S}\mathcal{M)}$ $=$ $e$-$Inj(_{S}%
\mathcal{M)}$ iff a semiring $S$ is a ring.\medskip

\noindent \textbf{Problem 2 }Is (1) of Theorem 4.5 true for all additively
regular semirings $S$?\medskip

\noindent \textbf{Problem 3 }Does an e-injective envelope exist for every
semimodule $M\in $ $|_{S}\mathcal{M}|$ over an additively idempotent
semiring $S$?\medskip

\noindent \textbf{Problem 4} Describe all semirings $S$ such that $Inj(_{S}%
\mathcal{M)}$ $\subset $ $e$-$Inj(_{S}\mathcal{M)}$.

\section{Characterizations of some special classes of semirings}

In \cite{il'in:dsoisadpopsos}, there were obtained semiring analogs of
well-known characterizations of (classical) semisimple, quasi-Frobenius, and
one-sided Noetherian rings by means of injective semimodules over them.
Motivated by this, it is quite natural to consider characterizations of
those classes of semirings in terms of e-injective semimodules over them.
And therefore, establishing the `e-injective' versions of those
characterizations constitutes the main goal of this section.\medskip

\noindent \textbf{Proposition 5.1 }(\textit{cf.} \cite[Theorem 1.2]%
{il'in:dsoisadpopsos})\textbf{\ }\textit{Let }$S$\textit{\ be a zerosumfree
semiring and }$\{M_{i}\}_{i\in I}$\textit{\ a family of left nonzero }$S$%
\textit{-semimodules. Then the left }$S$\textit{-semimodule }$\oplus _{i\in
I}M_{i}$\textit{\ is e-injective iff all semimodules }$\{M_{i}\}_{i\in I}$
\textit{are e-injective and }$|$\textit{\ }$I$\textit{\ }$|$\textit{\ }$%
<\infty $\textit{.\medskip }

\noindent \textit{Proof} \ $\Longrightarrow $. Let $M=\oplus _{i\in I}M_{i}$
be an e-injective left $S$-semimodule. By Corollary 3.3 (3), for every $i\in
I$, the $S$-semimodule $M_{i}$ is e-injective, and hence, by Proposition
4.1, there exists the element $z_{i}\in M_{i}$ such that $m_{i}+z_{i}=z_{i}$
for all $m_{i}\in M_{i}$. Then, the image of the canonical embedding $\mu
:M\rightarrowtail \overline{M}:=\Pi _{i\in I}$ $M_{i}$ is clearly a
subtractive subsemimodule, and therefore, there exists a homomorphism $%
\varphi :\overline{M}\longrightarrow M$ such that $\varphi \mu =1_{M}$.

Now let $z:=(z_{i})_{i\in I}\in \overline{M}$. Since $\varphi (z)\in M$,
there exists a finite subset $J\subseteq I$ and a collection of elements $%
(m_{j})_{j\in J}$ with $m_{j}\in M_{j}$, such that $\varphi (z)=\Sigma
_{j\in J}m_{j}$. We shall show that actually $J=I$. Indeed, for any index $%
i\in I\setminus J$, it is easy to see that $\mu (z_{i})+z=z$ and%
\begin{equation*}
\Sigma _{j\in J}m_{j}=\varphi (z)=\varphi (\mu (z_{i})+z)=\varphi (\mu
(z_{i}))+\varphi (z)=z_{i}+\Sigma _{j\in J}m_{k}\text{;}
\end{equation*}%
and as $i\notin J$, one has that $z_{i}=0$. Whence, for any $m_{i}\in M_{i},$
we have that $0=z_{i}=z_{i}+m_{i}=0+m_{i}=m_{i},$ \emph{i.e.}, $M_{i}=0$
what contradicts to $M_{i}\neq 0$. Thus, $J=I$.

$\Longleftarrow $. It follows straightforwardly from Corollary 3.3. \textit{%
\ \ \ \ \ \ }$_{\square }\medskip $

For given a semiring congruence $\theta $ on a semiring $S$ and a semimodule
$M\in $ $|_{S}\mathcal{M}|$, we say that an element $m\in M$ is \emph{%
compatible with} $\theta $ iff $s_{1}\theta s_{2}\Longrightarrow
s_{1}m=s_{2}m$ for all $s_{1},s_{2}\in S$; and let $M(\theta ):=\{m\in M\mid
m$ is compatible with $\theta \}$. As was shown in \cite[Propositions 1.2
and 1.3]{il'in:dsoisadpopsos}, $M(\theta )$ is a left $S/\theta $-semimodule
and $M(\theta )\simeq \mathrm{Hom}_{S}(S/\theta ,M)$ as left $S$%
-semimodules. From these observations and Corollary 3.5, we immediately
obtain the following e-injective analog of \cite[Proposition 1.4]%
{il'in:dsoisadpopsos}:\medskip

\noindent \textbf{Proposition 5.2 }\textit{If }$M$\textit{\ is an
e-injective left }$S$\textit{-semimodule, then }$M(\theta )$\textit{\ is an
e-injective left }$S/\theta $\textit{-semimodule.\medskip }

In the next result, we complement \cite[Theorem 3.4]{il'in:dsoisadpopsos}
and describe the class of semirings all of whose semimodules are
e-injective.\medskip

\noindent \textbf{Theorem 5.3 }\textit{The following conditions for a
semiring }$S$\textit{\ are equivalent:\medskip }

\textit{(1) All left }$S$\textit{-semimodules are e-injective;}

\textit{(2) }$S$\textit{\ is a classical semisimple ring.\medskip }

\noindent \textit{Proof} \ (1) $\Longrightarrow $ (2). By our hypothesis, $%
S^{(\mathbb{N})}$ is an e-injective left $S$-semimodule. Let $\theta $ be
the Bourne relation on $S$ corresponding to the ideal $V(S)$ (see 2.2 and
2.4). Then, it is easy to see that $S/\theta $ is a zerosumfree semiring,
and for the e-injective left $S$-semimodule $S^{(\mathbb{N})}$ there is the
natural $S/\theta $-semimodule isomorphism $S^{(\mathbb{N})}(\theta )\simeq
S(\theta )^{(\mathbb{N})}$, and therefore, by Proposition 5.2, $S^{(\mathbb{N%
})}(\theta )$ is an e-injective left $S/\theta $-semimodule. From the latter
and Proposition 5.1, $S(\theta )=0$.

Also, it is easy to verify that the image of the canonical embedding $\xi
_{V(S)}:V(S)\rightarrowtail D(S)$ (see 2.7) is a subtractive $S$%
-subsemimodule. For $_{S}S$ is e-injective, the natural injection $\mu
:V(S)\rightarrowtail S$ can be extended to an $S$-homomorphism $\psi
:D(S)\longrightarrow S$ such that $\psi \xi _{V(S)}=\mu $. Since $D(S)$ is a
ring, it is clear that $e:=\psi (1_{D(S)})\in V(S)$. For every $s\in V(S)$,
we have $s=\mu (s)=\psi (\xi _{V(S)}(s))=\psi (s1_{D(S)})=s\psi
(1_{D(S)})=se $. In particular, $e^{2}=e$, and so $%
s(1-e)=se(1-e)=s(e-e)=s0=0 $ for all $s\in V(S)$, and applying \cite[Lemma
1.1]{il'in:dsoisadpopsos}, we conclude that $1=e\in V(S)$, \emph{i.e.}, $S$
is actually a ring and, by \cite[Theorem 1.2.9]{lam:afcinr}, even a
semisimple ring as the concepts of \emph{injectivity} and \emph{e-injectivity%
} coincide for modules over rings.

(2) $\Longrightarrow $ (1). This follows immediately from \cite[Theorem 1.2.9%
]{lam:afcinr} as the concepts of \emph{injectivity} and \emph{e-injectivity}
coincide for modules over rings. \textit{\ \ \ \ \ \ }$_{\square }\medskip $

The celebrated Faith-Walker Theorem provides the characterization of
quasi-Frobenius rings as rings over which the classes of projective and
injective modules coincide (\emph{e.g.}, \cite[Theorem 31.9]{andfull:racom},
or \cite[Theorem 15.9]{lam:lomar}). In \cite[Theorem 3.5]%
{il'in:dsoisadpopsos}, this characterization has been generalized in the
semiring setting, and our next result is an `e-injective' version of \cite[%
Theorem 3.5]{il'in:dsoisadpopsos}.\medskip

\noindent \textbf{Theorem 5.4 }\textit{The following conditions for a
semiring }$S$\textit{\ are equivalent:\medskip }

\textit{(1) }$S$\textit{\ is a quasi-Frobenius ring;}

\textit{(2) All projective left }$S$\textit{-semimodules are e-injective;}

\textit{(3) All e-injective left }$S$\textit{-semimodules are projective,
and the }$S$\textit{-semimodule }$S/V(S)$\textit{\ can be embedded in an
e-injective left }$S$\textit{-semimodule.\medskip }

\noindent \textit{Proof} \ (1) $\Longrightarrow $ (2). This follows from
\cite[Theorem 31.9]{andfull:racom}, or \cite[Theorem 15.9]{lam:lomar}, and
the fact that the notions of \emph{injectivity} and \emph{e-injectivity}
coincide for modules over rings.

(2) $\Longrightarrow $ (1). Since the projective left $S$-semimodule $S^{(%
\mathbb{N})}$ is e-injective, as in the proof of the implication (1) $%
\Longrightarrow $ (2) of Theorem 5.4 we have that $S$ is a ring and then use
\cite[Theorem 31.9]{andfull:racom}.

(1) $\Longrightarrow $ (3). If $S$ is a ring, $S/V(S)=0$ is an e-injective $%
S $-semimodule. Therefore, the statement follows from \cite[Theorem 31.9]%
{andfull:racom}.

(3) $\Longrightarrow $ (1). Let $\theta $ be the Bourne relation on $S$
induced by $V(S)$. By the hypothesis, $S/\theta $ can be embedded in an
e-injective left $S$-semimodule $M$. It is easy to see that every element of
$S/\theta $ is compatible with $\theta $, and hence, $S/\theta $ can be
embedded in $M(\theta )$. Let $I$ be an infinite set such that $\left\vert
I\right\vert \geq \left\vert S\right\vert $. By Corollary 3.3 (2), $M^{I}$
is an e-injective left $S$-semimodule, and hence, a projective left $S$%
-semimodule; and therefore, by \cite[Theorem 2.1]{il'in:dsoisadpopsos}, $%
(M,+,0)$ is an abelian group. For $_{S}M$ is e-injective, by Proposition
5.2, $M(\theta )$ is an e-injective left semimodule over the zerosumfree
semiring $S/\theta $; and by Proposition 4.1, there exists an element $z\in
M(\theta )$ such that $z+m=z$ for all $m\in M(\theta )$. In particular, $%
z+z=z$, and, hence, as $M(\theta )\subseteq M$ and $(M,+,0)$ is a group, $%
z=0 $, and therefore, $M(\theta )=0$ and $S/\theta =0,$ \emph{i.e.} $S$ is a
ring. From the latter, the implication follows from \cite[Theorem 31.9]%
{andfull:racom}. \textit{\ \ \ \ \ \ }$_{\square }\medskip $

Also, the famous and very important characterization of Noetherian rings as
rings over which direct sum of injective modules is always an injective
module given by H. Bass and Z. Papp (see, \textit{e.g.}, \cite[Theorem 3.46]%
{lam:lomar}) has been generalized in the semiring setting in \cite[Theorem
3.6]{il'in:dsoisadpopsos}, and we conclude this section by an `e-injective'
version of the latter.\medskip

\noindent \textbf{Theorem 5.5 }\textit{The following conditions for a
semiring }$S$\textit{\ are equivalent:\medskip }

\textit{(1) }$S$\textit{\ is a left Noetherian ring;}

\textit{(2) Every direct sum of e-injective left }$S$\textit{-semimodules is
e-injective, and the }$S$\textit{-semimodule }$S/V(S)$\textit{\ can be
embedded in an e-injective left }$S$\textit{-semimodule.\medskip }

\noindent \textit{Proof} \ (1) $\Longrightarrow $ (2).\textit{\ }$S/V(S)$ is
obviously the zero e-injective left $S$-semimodule for any left Noetherian
ring $S$, and the implication follows from \cite[Theorem 3.46]{lam:lomar}.

(2) $\Longrightarrow $ (1). Let $\theta $ be the Bourne relation on $S$
induced by $V(S)$. By the hypothesis, $S/\theta $ can be embedded in an
e-injective left $S$-semimodule $M$. It is easy to see that every element of
$S/\theta $ is compatible with $\theta $, and hence, $S/\theta $ can be
embedded in $M(\theta )$. Then, as $M^{(\mathbb{N})}$ is an e-injective left
$S$-semimodule and applying Proposition 5.2, we have that $M(\theta )^{(%
\mathbb{N})}\simeq M^{(\mathbb{N})}(\theta )$ is an e-injective left $%
S/\theta $-semimodule. For $S/\theta $ is zerosumfree and applying
Proposition 5.1, we have $M(\theta )=0$ and, hence, $S/\theta =0$.
Therefore, by \cite[Theorem 3.46]{lam:lomar}, $S$ is a left Noetherian ring.
\textit{\ \ \ \ \ \ }$_{\square }\medskip $

\section{Semirings all of whose finitely generated semimodules are
e-injective}

As was shown in \cite[Corollary 3.2]{ikn:thstososaowcsap}, every semimodule
can be represented, in a canonical way, as a colimit of its cyclic
subsemimodules. This observation motivates studying of semirings over which
any semimodule is a colimit of cyclic semimodules possessing some special
properties (see, \emph{e.g.,} \cite{aikn:ovsasaowcsai}, \cite%
{ikn:thstososaowcsap}). Thus, it is quite natural that the main goal of this
section is to present complete characterizations of semirings, belonging to
some important special classes of semirings, all of whose cyclic (finitely
generated) semimodules are e-injective. But first we make some
interesting and useful general observations.\medskip

\noindent \textbf{Proposition 6.1 }\textit{The class of semirings all of
whose cyclic (finitely generated) left semimodules are e-injective is closed
under homomorphic images and finite products.\medskip }

\noindent \textit{Proof} \ Let $\mathfrak{C}$ ( $\mathfrak{F}$) be the class
of semirings all of whose cyclic (finitely generated) left semimodules are
e-injective. First consider the statement for the class $\mathfrak{C}$. So,
let $S,$ $T$ be semirings with $S\in \mathfrak{C}$ and $\pi
:S\longrightarrow T$ a surjective semiring homomorphism. We claim that $T\in
\mathfrak{C}$. By \cite[Section 4]{kat:thcos}, $\pi $ induces two functors
--- the \emph{restriction of scalars functor } $\pi ^{\#}:$ $_{T}\mathcal{M}%
\longrightarrow $ $_{S}\mathcal{M}$ and the \emph{extension functor} $\pi
_{\#}:=T\otimes _{S}-:$ $_{S}\mathcal{M}\longrightarrow $ $_{T}\mathcal{M}$
--- such that $\pi _{\#}\ $is a left adjoint to the functor $\pi ^{\#}$,
\textit{i.e.}, $\pi _{\#}\ \dashv $ $\pi ^{\#}$ \cite[Proposition 4.1]%
{kat:thcos}. Moreover, by \cite[Proposition 4.6]{kat:thcos}, $\pi _{\#}\pi
^{\#}\simeq Id_{_{T}\mathcal{M}}$; also, it is easy to see that a short
sequence $0\longrightarrow A\overset{f}{\longrightarrow }B\overset{g}{%
\longrightarrow }C\longrightarrow 0$ is exact in $_{T}\mathcal{M}$ iff it is
exact in $_{S}\mathcal{M}$, and any semimodule $M\in $ $|_{T}\mathcal{M}|$
is cyclic iff $M\in $ $|_{S}\mathcal{M}|$ is cyclic. By using these
observations and the assumption $S\in \mathfrak{C}$, it is easy to see that
the functor $\mathrm{Hom}_{T}(-,M):$ $_{T}\mathcal{M\longrightarrow }$ $%
\mathcal{M}$ preserves short exact sequences iff the functor $\mathrm{Hom}%
_{S}(-,M):$ $_{S}\mathcal{M\longrightarrow }$ $\mathcal{M}$ does it, and
therefore, $T\in \mathfrak{C}$, \textit{i.e.}, $\mathfrak{C}$ is closed
under homomorphic images. Using the same arguments, we obtain the closedness
of the class $\mathfrak{F}$ with respect of homomorphic images.

It is clear that it is enough to show only for two semirings $S_{1}$, $%
S_{2}\in \mathfrak{C}$ ( $\mathfrak{F}$) that $S:=S_{1}\oplus S_{2}\in
\mathfrak{C}$ ( $\mathfrak{F}$) as well. Obviously, any $S$-semimodule $M\in
$ $_{S}\mathcal{M}$ actually is the direct sum of its $S$-subsemimodules $%
S_{1}M:=\{(s_{1},0)$ $M\mid s_{1}\in S_{1}\}$ and $S_{2}M:=\{(0,s_{2})$ $%
M\mid s_{2}\in S_{2}\}$,\textit{\ i.e.}, $M=S_{1}M\oplus S_{2}M$. Moreover,
it is easy to see that any homomorphism $f:$ $M=S_{1}M\oplus S_{2}M$ $%
\mathcal{\longrightarrow }$ $S_{1}N\oplus S_{2}N=N$ between $S$-semimodules $%
M$, $N\in $ $_{S}\mathcal{M}$ is, in fact, the direct sum of the two
corresponding $S_{1}$-homomorphism $f_{1}:$ $S_{1}M$ $\mathcal{%
\longrightarrow }$ $S_{1}N$ and $S_{2}$-homomorphism $f_{2}:$ $S_{2}M$ $%
\mathcal{\longrightarrow }$ $S_{2}N$ between $S_{1}M$, $S_{1}N\in $ $_{S_{1}}%
\mathcal{M}$ and each $S_{2}M$, $S_{2}N\in $ $_{S_{2}}\mathcal{M}$,
respectively. From observations for any functor $\mathrm{Hom}_{S}(-,M):$ $%
_{S}\mathcal{M\longrightarrow }$ $\mathcal{M}$, we readily have $\mathrm{Hom}%
_{S}(-,M)=\mathrm{Hom}_{S}(-,S_{1}M\oplus S_{2}M)\simeq $ $\mathrm{Hom}%
_{S}(-,S_{1}M)\times \mathrm{Hom}_{S}(-,S_{2}M)\simeq $ $\mathrm{Hom}%
_{S_{1}}(-,S_{1}M)\times \mathrm{Hom}_{S_{2}}(-,S_{2}M)$, and therefore, by
Proposition 3.1, the functor $\mathrm{Hom}_{S}(-,M)$ is e-injective as soon
as the functors $\mathrm{Hom}_{S_{1}}(-,S_{1}M)$ and $\mathrm{Hom}%
_{S_{2}}(-,S_{2}M)$ are e-injective as well. Now it follows immediately that
$S\in \mathfrak{C}$ ( $\mathfrak{F}$) if we observe that the semimodules $%
S_{1}M\ $and $S_{2}M$ are cyclic, or finitely generated, provided that the
semimodule $M$ is itself cyclic, or finitely generated, respectively.
\textit{\ \ \ \ \ \ }$_{\square }\medskip $

\noindent \textbf{Proposition 6.2 }\textit{The following conditions for a
semiring }$S$\textit{\ are equivalent:\medskip }

\textit{(1) All cyclic (finitely generated) left }$S$\textit{-semimodules
are e-injective;}

\textit{(2) }$S\simeq R\oplus T$\textit{, where }$R$\textit{\ is a
(classical) semisimple ring, and }$T$\textit{\ is a zerosumfree semiring
with an infinite element whose all cyclic (finitely generated) left
semimodules are e-injective.\medskip }

\noindent \textit{Proof} \ (1) $\Longrightarrow $ (2). Let $\equiv _{V(S)}$
be the Bourne congruence on $S$. It is clear that the factor semiring $%
\overline{S}:=S/\equiv _{V(S)}$ is zerosumfree. So, by Proposition 6.1,
every cyclic (finitely generated) left $\overline{S}$-semimodule is
e-injective. In particular, $_{\overline{S}}\overline{S}$ is e-injective
and, by Proposition 4.1, contains an infinite element whence a zeroic
zerosumfree semiring. Whence, by \cite[Proposition 2.9]{il'in:v-s}, $%
S=R\oplus T$, where $R$ is a ring and $T$ is a semiring isomorphic to $%
\overline{S}$. Since $R$ is a homomorphic image of $S$, all cyclic (finitely
generated) left $R$-modules are injective by Proposition 6.1, and hence, by
the celebrated Osofsky's Theorem (\cite[Theorem, p. 649]{osof:raowfgmai}) $R$
is a semisimple ring.

(2) $\Longrightarrow $ (1). This follows immediately from Proposition 6.1
and \cite[Theorem; p. 649]{osof:raowfgmai}. \textit{\ \ \ \ \ \ }$_{\square
}\medskip $

As was shown in \cite{il'in:svwasai}, all finitely generated left $S$%
-semimodules over a semiring $S$ are injective if and only if a semiring $S$
is a (classical) semisimple ring. However, as we will demonstrate below, in
the general semiring setting the e-injective version of this result is not
true\emph{. }From Proposition 6.2, we see that the problem of describing
semirings all of whose cyclic (finitely generated) semimodules are e-
injective is\ actually reduced to the corresponding problem for the class of
zerosumfree semirings with infinite elements that, particularly, includes
the very important subclass of bounded distributive lattices. Therefore, it
is quite natural that a characterization of bounded distributive lattices
all of whose cyclic (finitely generated) semimodules are e-injective
constitutes our next goal, to achieve which we need first to establish the
following important facts.\medskip

\noindent \textbf{Lemma 6.3 }\textit{If }$S$\textit{\ is a bounded
distributive lattice all of whose cyclic semimodules are e-injective, then
it is a Boolean algebra.\medskip }

\noindent \textit{Proof} \ First let us show that for any $a\in S$, we have $%
Sa+Ann(a)=S$, where $Ann(a):=$ $(0:$ $a)$ is the \emph{annihilator} of the
element $a$. Indeed, suppose that $Sa+Ann(a)\varsubsetneq S$, then $%
Sa+Ann(a) $ is contained in a maximal ideal $I$ of $S$. Obviously, $I$, as
any ideal of any bounded distributive lattice, is strongly subtractive. It
implies that the set $\{0\}\cup \{sa\ |\ s\in S\setminus I\}$ is a submonoid
of the monoid $(S,+,0)$. Moreover, by \cite[Corollary 7.13]{golan:sata}, $I$
is also a prime ideal of $S$. If $sa\neq 0$ for all $s\in S\setminus I$,
then the monoid $\{0\}\cup \{sa\ |\ s\in S\setminus I\}$ becomes a left $S$%
-semimodule by setting $s(ra)=0$ for all $s\in I$ and $r\notin I$, and $%
s(ra)=sra$ for all $s\notin I$ and $r\notin I$. Next let us take an
arbitrary element $m\notin S$ and on the set $V:=\{0\}\cup \{sa\ |\ s\in
S\setminus I\}\cup \{m\}$ extend the operations of the $S$-semimodule $%
\{0\}\cup \{sa\ |\ s\in S\setminus I\}$ by setting $m+0=0+m=m+m=m$, $%
m+sa=sa+m=sa$ for all $s\in S\setminus I$, and $sm=0$ for all $s\in I$, and $%
sm=m$ for all $s\in S\setminus I$. One can easily verify that $V$ becomes a
left $S$-semimodule. Obviously, $M:=\{0,m\}$ is a subtractive $S$%
-subsemimodule of $V$. Clearly, $M$ is cyclic and, hence, e-injective;
therefore, there exists an $S$-homomorphism $f:V\longrightarrow M$ such that
$f|_{M}=id_{M}$. It implies that $f(sa)=m$ for all $s\in S\setminus I$. In
particular, we have that $f(a)=m$. On the other hand, since $%
Sa+Ann(a)\subseteq I$, we get that $a\in I$, and hence, $%
f(a)=f(a^{2})=af(a)=am=0$. Whence, $sa=0$ for some $s\in S\setminus I$, that
is, $s\in Ann(a)\subseteq I$, and hence, $s\in I$; and therefore, this
contradiction implies $S=Sa+Ann(a)$ for all $a\in S$.

Now we show that $S$ is a Boolean algebra. Indeed, for any $a\in S$, we have
$S=Sa+Ann(a)$; and hence, $1=sa+x$ for some $s\in S$ and $x\in Ann(a)$. We
then have that
\begin{equation*}
1=a+1=a+sa+x=a(1+s)+x=a1+x=a+x.
\end{equation*}%
For $x\in Ann(a)$, we have $xa=0$, and therefore, $x$ is the complement to $%
a $, and $S$ is a Boolean algebra. \textit{\ \ \ \ \ \ }$_{\square }\medskip
$

Moreover, for the semiring $S$ in the previous lemma, we can make a more
precise observation, namely:\medskip

\noindent \textbf{Lemma 6.4 }\textit{If }$S$\textit{\ is a bounded
distributive lattice all of whose cyclic semimodules are e-injective, then }$%
S$\textit{\ is a complete Boolean algebra.\medskip }

\noindent \textit{Proof} \ By \cite[Theorem X.9]{birkhoff:lathe} (or, \cite[%
Theorem 8.5]{skornyakov:ets}) and Lemma 6.3, the Boolean algebra $S$ can be
considered as a Boolean subalgebra of a complete Boolean algebra $B$; and
hence, $S\subseteq B$, $0_{S}=0_{B}$, and $1_{S}=1_{B}$, and $B$ is
naturally a left $S$-semimodule. Define on the left $S$-semimodule $S\times
B $ the relation $\sim $ as follows: $(s_{1},b_{1})\sim (s_{2},b_{2})$ iff
1) $b_{1}=b_{2}$ and 2)~there exist $n\in \mathbb{N}$, $x_{i},x_{i}^{\prime
}\in B$, $s_{i}^{\prime },s_{i}^{\prime \prime }\in S$, $i=1,\ldots ,n$,
such that $x_{i}+x_{i}^{\prime }=b_{1}$ for all $i$, and $%
s_{1}+\sum_{i=1}^{n}s_{i}^{\prime }x_{i}=s_{2}+\sum_{i=1}^{n}s_{i}^{\prime
\prime }x_{i}$. Repeating verbatim the proof of \cite[Proposition 1.6]%
{il'in:dsoisadpopsos}, one easily sees that $\sim $ is a congruence on $%
_{S}(S\times B)$.

Let $\widetilde{S}:=(S\times B)/{\sim }$ and $[s,b]$ be the class containing
the pair $(s,b)$. It is easy to see that the image of the embedding $\mu
:S\rightarrowtail \widetilde{S}$, $s\longmapsto \lbrack s,0]$, is a
subtractive $S$-semimodule of $_{S}\widetilde{S}$. Hence, for $_{S}S$ is
e-injective, there exists an $S$-homomorphism $\varphi \colon \widetilde{S}%
\rightarrow S$ such that $\varphi \mu =1_{S}$. By the dual of \cite[Theorem
IX.2]{birkhoff:lathe}, in order to prove $S$ is complete, it is sufficient
to show that every family $\{s_{i}\}_{i\in I}$ of elements in $S$ has a
least upper bound $s\in S$. Since $S\subseteq B$ and $B$ is complete, a
family $\{s_{i}\}_{i\in I}$ has a least upper bound $b\in B$ in $B$, and let
$s:=\varphi ([0,b])$. For each $i\in I$, it is easy to see that $%
[s_{i},b]=[0,b]$: Indeed, $(s_{i},b)\sim (0,b)$ by putting $n=1$, $x_{1}=b$,
$x_{1}^{\prime }=0$, $s_{1}^{\prime }=0$, $s_{1}^{\prime \prime }=s_{i}$.
Therefore,
\begin{equation*}
s_{i}+s=\varphi (\mu (s_{i}))+\varphi ([0,b])=\varphi
([s_{i},0]+[0,b])=\varphi ([s_{i},b])=\varphi ([0,b])=s\text{,}
\end{equation*}%
and hence, $s$ is an upper bound for $\{s_{i}\}_{i\in I}$ in $S$. If $y\in S$
is another upper bound for $\{s_{i}\}_{i\in I}$, then $b\leq y$ as $b$ is
the least upper bound for $\{s_{i}\}_{i\in I}$ in $B$; and therefore, $by=b$
and $sy=\varphi ([0,b])y=\varphi ([0,by])=\varphi \lbrack 0,b]=s$, \textit{%
i.e.}, $s\leq y$ and $s$ is a least upper bound for $\{s_{i}\}_{i\in I}$ in $%
S$. \textit{\ \ \ \ \ \ }$_{\square }\medskip $

Now we are ready to characterize bounded distributive lattices all of whose
cyclic (finitely generated) semimodules are e-injective, namely:\medskip

\noindent \textbf{Theorem 6.5 }\textit{The following conditions for a
bounded distributive lattice }$S$\textit{\ are equivalent:\medskip }

\textit{(1) All finitely generated }$S$\textit{-semimodules are e-injective;}

\textit{(2) All cyclic }$S$\textit{-semimodules are e-injective;}

\textit{(3) \medskip }$S$\textit{\ is a finite Boolean algebra.}

\noindent \textit{Proof} (1) $\Longrightarrow $ (2). It is obvious.

(2) $\Longrightarrow $ (3). By Lemma 6.4, $S$ is a complete Boolean algebra.
Suppose that $S$ is infinite, then, by \cite[Lemma 4.2]{aikn:ovsasaowcsai}
for example, it contains a countable set of orthogonal idempotents $%
\{e_{n}\,|\,n\in \mathbb{N}\}$, and we have the ideal $I:=\Sigma _{n\in
\mathbb{N}}$ $Se_{n}$ and the factor algebra $\overline{S}:=S/I$. For $%
\overline{S}$ is a Boolean algebra and by Proposition 6.1, we have that all
cyclic $\overline{S}$-semimodules are e-injective and, by Lemma 6.4 again, $%
\overline{S}$ is a complete Boolean algebra. On the other hand, repeating
verbatim the proof of \cite[Theorem 4.3]{aikn:ovsasaowcsai}, we have that $%
\overline{S}$ is not a complete Boolean algebra, and from this contradiction
we conclude that $S$ is a finite Boolean algebra.

(3) $\Longrightarrow $ (1). Let $S$ be a finite Boolean algebra and $%
\{e_{1},e_{2},...,e_{n}\}$ the set of all \emph{atoms} of $S$. Clearly, $%
S=Se_{1}\oplus Se_{2}\oplus ...\oplus Se_{n}$ and $Se_{i}=\{0,e_{i}\}\simeq
\mathbf{B}$ for each $i$. Then, applying Theorem 4.4 (5) and Proposition
6.1, we conclude the proof. \textit{\ \ \ \ \ \ }$_{\square }\medskip $

As a corollary of Theorem 6.5, we are able to extend the characterization of
finite Boolean algebras among bounded distributive lattices given in \cite[%
Corollary 4.4]{aikn:ovsasaowcsai}:\medskip

\noindent \textbf{Corollary 6.6} \textit{The following conditions for a
bounded distributive lattice }$S$\textit{\ are equivalent:\medskip }

\textit{(1) All cyclic }$S$\textit{-semimodules are projective;}

\textit{(2) All subsemimodules of the regular semimodule }$_{S}S$\textit{\
are injective;}

\textit{(3) All cyclic }$S$\textit{-semimodules are e-injective;}

\textit{(4)} \textit{All cyclic }$S$\textit{-semimodules are injective;}

\textit{(5) }$S$\textit{\ is a finite Boolean algebra.\medskip }

As another consequence of Theorem 6.5, we obtain a complete description of
subtractive semirings all of whose finitely generated semimodules are
e-injective.\medskip

\noindent \textbf{Theorem 6.7 }\textit{The following conditions for a left
subtractive semiring }$S$\textit{\ are equivalent:\medskip }

\textit{(1) All finitely generated }$S$\textit{-semimodules are e-injective;}

\textit{(2) All cyclic }$S$\textit{-semimodules are e-injective;}

\textit{(3) \medskip }$S\simeq R\oplus T$\textit{, where }$R$\textit{\ is a
semisimple ring and }$T$\textit{\ is a finite Boolean algebra.}

\noindent \textit{Proof} (1) $\Longrightarrow $ (2). It is obvious.

(2) $\Longrightarrow $ (3). By Proposition 6.2, $S\simeq R\oplus T$, where $%
R $ is a semisimple ring and $T$ a semiring with an infinite element $\infty
$ such that all cyclic left $T$-semimodules are e-injective. By \cite[Lemma
4.7]{knt:ossss}, $T$ is a left subtractive semiring. Whence, the left ideal $%
T\infty $ is a subtractive ideal of $T$ and it follows from $1_{T}+\infty
=\infty $ that $1_{T}\in T\infty $, \textit{i.e.}, $t\infty =1_{T}$ for some
$t\in T$. From the latter and $\infty ^{2}+\infty =\infty $, we have $\infty
=\infty +1_{T}=1_{T}$, and hence, $x+1_{T}=1_{T}$ for all $x\in T$.

As, for each $a\in T$, the cyclic left $T$-semimodule $Ta$ is both
e-injective and subtractive subsemimodule of $T$ by our hypothesis.
Therefore, there exists a homomorphism $\varphi :T\longrightarrow Ta$ such
that $\varphi |_{Ta}=1_{Ta}.$ In particular, $a=\varphi (a)=\varphi
(a1)=a\varphi (1)$ and $\varphi (1)\in Ta$. Hence, $\varphi (1)=xa$ for some
$x\in T$ and $a=axa$. From the latter, repeating verbatim the proof of \cite[%
Theorem 4.6]{aikn:ovsasaowcsai}, we immediately obtain that $T$ is a bounded
distributive lattice, and therefore, applying Theorem 6.5, we conclude that $%
T$ is a finite Boolean algebra.

(3) $\Longrightarrow $ (1). This follows immediately from Proposition 6.1,
Theorem 6.5, and \cite[Theorem, p. 649]{osof:raowfgmai}. \textit{\ \ \ \ \ \
}$_{\square }\medskip $

By \cite[Corollary 3.2]{ikn:thstososaowcsap}, every $S$-semimodule over a
semiring $S$ can be canonically represented as a colimit of a diagram of
cyclic $S$-semimodules. In light of this observation, Corollary 6.6 and
Theorem 6.7, it seems to be reasonable and interesting the following
problem.\medskip

\noindent \textbf{Problem 5} Describe all semirings $S$ for which the
classes of cyclic injective and cyclic e-injective semimodules
coincide.\medskip

Now, in what follows, we use the concepts and notations from 2.4, 2.5 and
2.6, and following \cite{Abu2014-b}, we say that a right $S$-semimodule $%
P_{S}$ is \textit{e-flat} iff the functor $P\otimes _{S}-$ $:$ $_{S}\mathcal{%
M}\longrightarrow \mathcal{M}$ preserves short exact sequences or
equivalently iff the functor $P\otimes _{S}-$ $\ $preserves subtractive
subsemimodules.\medskip\

\noindent \textbf{Lemma 6.8 }\textit{Every projective right }$S$\textit{%
-semimodule is e-flat.\medskip }

\noindent \textit{Proof} \ As was shown in \cite{kat:tpaieosoars} and \cite%
{kat:thcos}, any tensor product functor $P\otimes _{S}-:$ $_{S}\mathcal{M}%
\longrightarrow \mathcal{M}$ has a right adjoint, and therefore, by \cite[%
The dual Theorem 5.5.1]{macl:cwm}, preserves colimits, in particular,
coproducts. Then, for projective semimodules are retracts of free ones, the
statement follows right away from \cite[Proposition 3.8]{kat:tpaieosoars}
and Corollary 3.3. \textit{\ \ \ \ \ \ }$_{\square }\medskip $

\noindent \textbf{Lemma 6.9 }\textit{Let }$F:$\textit{\ }$_{S}\mathcal{M}%
\rightleftarrows $\textit{\ }$_{T}\mathcal{M}:G$\textit{\ be an equivalence
between the semimodule categories }$_{S}\mathcal{M}$ and $_{T}\mathcal{M}$%
\textit{. Then a left }$S$\textit{-semimodule }$_{S}M\in $\textit{\ }$_{S}%
\mathcal{M}$ \textit{is e-injective iff }$_{T}F(M)\in $\textit{\ }$_{T}%
\mathcal{M}$\textit{\ is e-injective.\medskip }

\noindent \textit{Proof} \ By \cite[Theorems 4.5 and 4.12]{kn:meahcos}, $%
F\simeq P\otimes _{S}-$ for some $(T,S)$-bisemimodule $P$ such that $P_{S}$
is a progenerator, $P^{\ast }:=\mathrm{Hom}_{S}(P_{S},S_{S})$ is a
progenerator in $\mathcal{M}_{T}$ and $T\simeq End(P_{S})$ as semirings, and
$G\simeq P^{\ast }\otimes _{T}-$. Since projective semimodules are e-flat by
Lemma 6.8, both $F$ and $G$ preserve short exact sequences, and the
statement readily follows from the natural functor isomorphisms $FG\simeq
Id_{_{S}\mathcal{M}}$ and $GF\simeq Id_{_{T}\mathcal{M}}$ \cite[Sect. 4.4]%
{macl:cwm} (see also the dual of \cite[Lemma 4.10]{kn:meahcos}). \textit{\ \
\ \ \ \ }$_{\square }\medskip $

We conclude these section and paper by giving a complete characterization of
simple semirings all of whose finitely generated left semimodules are
e-injective.\medskip

\noindent \textbf{Theorem 6.10 }\textit{The following conditions for a
semiring }$S$\textit{\ are equivalent:\medskip }

\textit{(1) }$S$\textit{\ is a simple semiring all of whose finitely
generated left semimodules are e-injective;}

\textit{(2) }$S$\textit{\ is a simple semiring all of whose cyclic left
semimodules are e-injective;}

\textit{(3) }$S$\textit{\ is isomorphic either to a matrix semiring }$%
M_{n}(D)$\textit{\ for some division ring }$D$\textit{\ and }$n\geq 1$%
\textit{, or to an endomorphism semiring }$End(L)$\textit{\ of a nonzero
finite distributive lattice }$L$\textit{.}$\medskip $

\noindent \textit{Proof} (1) $\Longrightarrow $ (2). It is obvious.

(2) $\Longrightarrow $ (3). By Proposition 6.2, $S$ is a simple semisimple
ring, or $S$ is a simple semiring with an infinite element. If $S$ is a
simple semisimple ring, then $S\simeq M_{n}(D)$ for some division ring $D$
and $n\geq 1$ by the classical Wedderburn-Artin structure theorem for rings
(see, \textit{e.g.}\emph{,} \cite[Wedderburn-Artin Theorem 3.5]{lam:afcinr}%
). Otherwise, $S$ is a simple semiring with an infinite element $\infty .$
In particular, $\infty +\infty =\infty $, whence $\infty \in I^{+}(S)$ and
hence $I^{+}(S)$ is a nonzero ideal of $S$. For $S$ is a simple semiring, $%
I^{+}(S)=S$ and $S$ is an additively idempotent simple semiring containing
an infinite element, and therefore, by \cite[Theorem 5.7]{knz:ososacs}, $%
S\simeq End(L)$ for some nonzero finite distributive lattice $L$.

(3) $\Longrightarrow $ (1). \textbf{Case I:} Let $S\simeq M_{n}(D)$ for some
division ring $D$ and $n\geq 1$. It follows by the Wedderburn-Artin Theorem
and the celebrated Osofsky's result \cite[Theorem]{osof:raowfgmai}, all
finitely generated left $S$-semimodules are (e-)injective.

\textbf{Case II: }Let\textbf{\ }$S\simeq End(L)$ for some nonzero finite
distributive lattice $L$. Then, by \cite[Theorem 5.7]{knz:ososacs}, $S$ is a
simple semiring that Morita equivalent to the Boolean semiring $\mathbf{B}$;
and let the functors $F:$ $_{S}\mathcal{M}\rightleftarrows $ $_{\mathbf{B}}%
\mathcal{M}:G$ establish an equivalence between the semimodule categories $%
_{S}\mathcal{M}$ and $_{\mathbf{B}}\mathcal{M}$. For any finitely generated
left $S$-semimodule $M\in $ $|_{S}\mathcal{M}|$, by \cite[Proposition 4.8]%
{kn:meahcos} and Theorem 4.4 (5), the semimodule $F(M)\in $ $|_{\mathbf{B}}%
\mathcal{M}|$ is a finitely generated e-injective $\mathbf{B}$-semimodule.
Then, applying Lemma 6.9 and the natural isomorphism $M\simeq G(F(M))$, we
have that the semimodule $M\in $ $|_{S}\mathcal{M}|$ is e-injective as well
and end the proof. \textit{\ \ \ \ \ \ }$_{\square }\medskip $

\textbf{Acknowledgement:} The first and fourth authors would like to
acknowledge a support provided by the Deanship of Scientific Research (DSR)
at King Fahd University of Petroleum $\&$ Minerals (KFUPM) for funding this
work through project RG1304.


\begin{thebibliography}{99}
\bibitem{Abu2014-a} {\small J. Abuhlail, Exact sequences of commutative
monoids and semimodules, \emph{Homology Homotopy Appl.} 16 (1) (2014),
199-214.}

\bibitem{Abu2014-b} {\small J. Abuhlail, Some remarks on tensor products and
flatness of semimodules, \emph{Semigroup Forum} 88 (3) (2014), 732-738.}

\bibitem{aikn:ovsasaowcsai} {\small J. Y. Abuhlail, S. N. Il'in, Y. Katsov,
T. G. Nam, On $V$-Semirings and Semirings all of whose Cyclic Semimodules
are Injective, \emph{Communications in Algebra}, 43 (2015), 4632--4654.}

\bibitem{andfull:racom} {\small F. W. Anderson and K. R. Fuller, \textit{%
Rings and Categories of Modules}, 2nd ed., Springer-Verlag, New York-Berlin,
1979. }

\bibitem{bshhurtjankepka:scs} {\small R. El Bashir, J. Hurt, A. Jan\v{c}a%
\v{r}\'{\i}k, and T. Kepka, Simple Commutative Semirings, \textit{J. Algebra}%
, \textbf{236 }(2001), 277--306. }

{\small
}

\bibitem{birkhoff:lathe} {\small G. Birkhoff, \textit{Lattice Theory},
American Mathematical Society, Providence, 1967. }

{\small
}

{\small
}

\bibitem{Bor1994a} F. Borceux, Handbook of Categorical Algebra. I, Basic
Category Theory, \emph{Cambridge Univ. Press} (1994).

\bibitem{cc:sofazf} {\small A. Connes and C. Consani, Schemes over $\mathbb{F%
}_1$ and zeta functions, \textit{Compos. Math.}, \textbf{146} (2010),
1383--1415. }

\bibitem{gg:eotv} {\small J. Giansiracusa and N. Giansiracusa, Equations of
tropical varieties, to appear in Duke Mathematical Journal (see also
preprint: arXiv:1308.0042v2). }

\bibitem{fofan:ipoba} {\small T. S. Fofanova, Injectivity of polygons over
Boolean algebras, \textit{Siberian Math. J.}, \textbf{13} (1972), 452--458
(in Russian). }

\bibitem{fof:podl} {\small T. S. Fofanova, Polygons over distributive
lattices, in: \textit{Universal Algebra, Colloq. Math. Soc. J\'{a}nos Bolyai
\# 29, }North-Holland Publishing Co., Amsterdam, 1982, 289--292. }

{\small
}

\bibitem{gelman:moha} S. I. Gelfand, Y. I. Manin, \textit{Methods of
homological algebra}, Second edition, Springer Monographs in Mathematics.
Springer-Verlag, Berlin, 2003

\bibitem{glaz:agtl} K. G\l azek(2001). \textit{A Guide to the Literature on
Semirings and their Applications in Mathematics and Information Science,}
Kluwer Academic Publishers, Dordrecht-Boston-London, 2001

\bibitem{golan:sata} {\small J. S. Golan, \textit{Semirings and their
Applications}, Kluwer Academic Publishers, Dordrecht-Boston-London, 1999. }

{\small
}

{\small
}

\bibitem{hk:tcos} {\small A. Horn, N. Kimura, The category of semilattices,
\textit{Algebra Universalis}, \textbf{1} (1971), 26--38. }

\bibitem{il'in:svwasai} {\small S. N. Il'in, Semirings over Which all
Semimodules are Injective (Projective), \textit{Matem. Vestn. Pedvuzov i
Univ. Volgo-Vyatsk. Regiona}, \textbf{8} (2006), 50--53. }

{\small
}

\bibitem{il'in:dsoisadpopsos} {\small S. N. Il'in, Direct Sums of Injective
Semimodules and Direct Products of Projective Semimodules Over Semirings,
\textit{Russian Mathematics}, \textbf{54} (2010), 27 - 37. }

\bibitem{il'in:v-s} {\small S. N. Il'in, V-semirings, \textit{Siberian
Mathematical Journal}, \textbf{53} (2012), 222 - 231. }

\bibitem{ikn:thstososaowcsap} {\small S. N. Il'in, Y. Katsov, T.G. Nam,
Toward Homological Structure Theory of Semimodules: On Semirings All of
Whose Cyclic Semimodules Are Projective, preprint: arXiv:1509.02997v1. }

\bibitem{ir:sa} {\small Z. Izhakian, L. Rowen, Supertropical algebra, \emph{%
Adv. Math.}, \textbf{225}(4) (2010), 2222--2286. }

\bibitem{jun:ccoss} {\small J. Jun, \u{C}ech cohomology of semiring schemes,
preprint: arXiv:1503.01389v1. }

\bibitem{kat:tpaieosoars} {\small Y. Katsov, Tensor products and injective
envelopes of semimodules over additively regular semirings, \textit{Algebra
Colloquium}, \textbf{4} (1997), 121-131. }

\bibitem{kat:ofsos} {\small Y. Katsov, On flat semimodules over semirings,
\textit{Algebra Universalis}, \textbf{51} (2004), 287-299. }

\bibitem{kat:thcos} {\small Y. Katsov, Toward homological characterization
of semirings: Serre's conjecture and Bass's perfectness in a semiring
context, \textit{Algebra Universalis}, \textbf{52} (2004), 197-214. }

\bibitem{knt:ossss} {\small Y. Katsov, T. G. Nam, N. X. Tuyen, On
subtractive semisimple semirings, \textit{Algebra Colloquium}, \textbf{16}
(2009), 415-426. }

\bibitem{kn:meahcos} {\small Y. Katsov and T. G. Nam, Morita Equivalence and
Homological Characterization of Semirings, \textit{J. Algebra Appl.},
\textbf{10} (2011), 445 - 473.}

\bibitem{knt:mosssparp} {\small Y. Katsov, T. G. Nam, N. X. Tuyen, More on
Subtractive Semirings: Simpleness, Perfectness and Related Problems, \textit{%
Comm. Algebra}, \textbf{39} (2011), 4342 - 4356. }

\bibitem{knz:ososacs} {\small Y. Katsov, T. G. Nam, J. Zumbr\"{a}gel, On
Simpleness of Semirings and Complete Semirings, \emph{J. Algebra Appl.},
\textbf{13}: 6 (2014). DOI: 10.1142/S0219498814500157. }

\bibitem{kn:orosarp} Y. Katsov, T. G. Nam, On radicals of semirings and
related problems, \textit{Comm. Algebra} \textbf{42} (2014), no. 12,
5065--5099.

\bibitem{knz:solpawcics} {\small Y. Katsov, T. G. Nam, J. Zumbr\"{a}gel, }%
Simpleness of Leavitt path algebras with coefficients in commutative
semiring, \textit{Semigroup Forum}, 2016. DOI 10.1007/s00233-016-9781-1.

\bibitem{kz:fsais} {\small A. Kendziorra and J. Zumbr\"{a}gel, Finite simple
additively idempotent semirings, \emph{J. Algebra}, \textbf{388} (2013),
43--64. }

\bibitem{kkn:srwz} {\small T. Kepka, J. Kortelainen and P. N\u{e}mec, Simple
semirings with zero, \emph{J. Algebra Appl.}. \textbf{15}: 3 (2016) 1650047. DOI:
10.1142/S021949881650047X. }

\bibitem{kn:srwlmae} {\small T. Kepka and P. N\u{e}mec, Simple semirings
with left muliplicatively absorbing elements, \emph{Semigroup Forum},
\textbf{91} (2015), 159--170. }

\bibitem{kilp-kn-mik:maac} {\small M. Kilp, U. Knauer, A. V. Mikhalev,
\textit{Monoids, Acts and Categories,} Walter de Gruyter, Berlin-New York,
2000. }

\bibitem{lam:lomar} {\small T. Y. Lam, \textit{Lectures on Modules and Rings,%
} Springer-Verlag, New York-Berlin, 1999. }

\bibitem{lam:afcinr} {\small T. Y. Lam, \textit{A first course in
noncommutative rings, }2nd Ed., Springer-Verlag, New York-Berlin, 2001. }

\bibitem{lor:tgob} {\small O. Lorscheid, The geometry of blueprints. Part I:
Algebraic background and scheme theory, \emph{Adv. Math.}, \textbf{229}(3)
(2012), 1804--1846. }

\bibitem{macl:cwm} {\small S. Mac Lane, \textit{Categories for the Working
Mathematician}, Springer-Verlag, New York-Berlin, 1971. }

\bibitem{mmr:pkcbosa} {\small G. Maze, C. Monico, and J. Rosenthal, Public
key cryptography based on semigroup actions, \emph{Adv. Math. Commun.},
\textbf{1} (2007), 489--507.
}

\bibitem{monico:ofcss} {\small C. Monico, On finite congruence-simple
semirings, \textit{J. Algebra}, \textbf{271} (2004), 846--854. }

{\small 
}

\bibitem{osof:raowfgmai} {\small B. L. Osofsky, Rings all of whose finitely
generated modules are injective, \emph{Pacific J. Math.}, \textbf{14}
(1964), 645 - 650. }

\bibitem{p:eosatfe} {\small A. Patchkoria, Extensions of Semimodules and
Takahashi functor $Ext_{\Lambda }(C,A)$, \emph{Homology, Homotopy and Appl.,}
\textbf{5} (2013), 387 - 406. }

\bibitem{rst:fsitg} {\small J. Richter-Gebert, B. Sturmfels, T. Theobald,
\emph{First steps in tropical geometry, in: Idempotent Mathematics and
Mathematical Physics}, in: Contemp. Math., vol. 377, American Mathematical
Society, Providence, RI, 2005, pp. 289--317. }

\bibitem{schubert:c} H. Schubert, \textit{Categories}, Springer-Verlag, New
York-Heidelberg, 1972.

\bibitem{skornyakov:ets} L. A. Skornyakov, L. A. \textit{\`{E}lementy teorii
struktur }(Russian) [\textit{Elements of lattice theory}], Second edition,
\textquotedblleft Nauka\textquotedblright , Moscow, 1982.

\bibitem{Tak81} {\small M. Takahashi, On the bordism categories. II.
Elementary properties of semimodules, \textit{Kobe J. Math.}, \textbf{9}
(1981), 495-530. }

\bibitem{Tak1982a} M. Takahashi, On the bordism categories. III. Functors
Hom and for semimodules, \emph{Math. Sem. Notes Kobe Univ.} 10 (2) (1982),
551-562.

\bibitem{Tak82} {\small M. Takahashi, Completeness and c-cocompleteness of
the category of semimodules, \textit{Kobe J. Math.}, \textbf{10} (1982),
551--562. }

\bibitem{Tak:eos} {\small M. Takahashi, Extensions of semimodules I, \textit{%
Kobe J. Math.}, \textbf{10} (1982), 563--592. }

{\small
}

\bibitem{zumbr:cofcsswz} {\small J. Zumbr\"{a}gel, Classification of finite
congruence-simple semirings with zero, \textit{J. Algebra Appl.}, \textbf{7}
(2008), 363--377. }

{\small
}

{\small
}

{\small 
}
\end{thebibliography}
\end{document}